\documentclass[12pt,a4paper, twoside]{article}
\usepackage{amsmath,amssymb,amsthm,amsfonts,mathrsfs,amscd,environ}
\usepackage{dsfont}
\usepackage{latexsym,enumerate,color,geometry,extarrows}

\usepackage{verbatim,fancyhdr,enumitem}%for comment paragraph fancyhdr,
 \NewEnviron{ews}{%
\begin{equation}\begin{split}
  \BODY
\end{split}\end{equation}
}

\NewEnviron{ews*}{%
\begin{equation*}\begin{split}
  \BODY
\end{split}\end{equation*}
}

\newtheorem{prop}{Proposition}[section]
 \newtheorem{thm}{Theorem}[section]
 \newtheorem{lem}[thm]{Lemma}

 \newtheorem{prp}[thm]{Proposition}
 \newtheorem{defn}{Definition}[section]
 \newtheorem{rem}{Remark}[section]
 
 \numberwithin{equation}{section}

\def\dif{{\mathord{{\rm d}}}}
\def\no{\nonumber}

\def\mR{{\mathbb R}}
\def\mE{{\mathbb E}}

\def\mA{{\mathbb A}}

\def\mD{{\mathbb D}}
\def\mE{{\mathbb E}}
\def\mF{{\mathbb F}}

\def\mH{{\mathbb H}}

\def\mR{{\mathbb R}}
\def\mS{{\mathbb S}}

\def\sF{{\mathscr F}}

\def\sB{{\mathscr B}}

\def\sF{{\mathscr F}}

\def\sL{{\mathscr L}}

\def\sP{{\mathscr P}}

\def\sS{{\mathscr S}}

\def\bd{\begin{defn}}
\def\ed{\end{defn}}
\def\bp{\begin{prop}}
\def\ep{\end{prop}}
\def\bl{\begin{lem}}
\def\el{\end{lem}}
\def\bt{\begin{thm}}
\def\et{\end{thm}}
\def\br{\begin{rem}}
\def\er{\end{rem}}
\allowdisplaybreaks
\title{{\bf Large deviation for slow-fast McKean-Vlasov stochastic differential equations driven by fractional Brownian motions and Brownian motions}
}
\author{
{\bf  Hao Wu$^{a)}$, Junhao Hu$^{a)}$, Chenggui Yuan$^{b)}$
 }\\
\footnotesize{$^{a)}$School of Mathematics and Statistics,
South-Central University For Nationalities}\\
\footnotesize{ Wuhan, Hubei 430000, P.R.China}\\
\footnotesize{Email: wuhaomoonsky@163.com},
\footnotesize{ junhaohu74@163.com}\\
\footnotesize{$^{b)}$Department of Mathematics, Swansea University, Bay campus, SA1 8EN, UK}\\
\footnotesize{Email: C.Yuan@Swansea.ac.uk}\\
}

\begin{document}

\allowdisplaybreaks
%\title{\huge\CJKfamily{song}Large deviation principle for neutral type stochastic differential equation with infinite memory
%\author{***
% }
%}
\maketitle
\begin{abstract}
In this article, we consider slow-fast McKean-Vlasov stochastic differential equations driven by fractional Brownian motions and  Brownian motions. We give a definition of the large deviation principle (LDP)  on the product space related to fractional Brownian motion and  Brownian motion, which is different from the traditional definition for LDP. Under some proper assumptions on coefficients,  LDP is investigated    for  this type of equations by  using the weak convergence method.
\end{abstract}\noindent

\noindent


\section{Introduction}   The LDP theory is used to solve the  asymptotic behaviour of rare events and  has a wide range of applications  such as in finance, statistic mechanics, biology, etc., see, e.g., \cite{DGC, T,PS,P}  and references therein. From the literature, we can see that there exist two main methods to study the LDP, one is based on contraction principle, that is, it relies on approximation arguments and exponential-type probability estimates,  see, e.g., \cite{F,LZ,RZ} and references therein.  The other one is the weak convergence method  due to \cite{DE}. This method has been proved to be very effective for slow-fast  stochastic models, etc., see, e.g., \cite{WRD,GSS,GG,PK,YXJ,GL}.% and references therein.

On the other hand, many researchers are interested in   Mckean-Vlasov stochastic differential equations (SDEs) driven by Brownian motion , in which the coefficients   depend not only on  the state process but also on its distribution.  McKean-Vlasov SDEs, being clearly more involved than It\^{o}'s SDEs, arise in
 McKean \cite {M2}, who was inspired by Kac's Programme in Kinetic Theory \cite{K}, as well as in some other areas of high interest such as
propagation of chaos phenomenon, PDEs, stability, invariant probability measures, social science, economics, engineering, etc. (see  e.g. \cite{CD, DQ2, C, MH, MMLM, Z1, GA,S}).  Later, the averaging principle and LDP for slow-fast McKean-Vlasov SDEs driven by  Brownian motions have been studied in many papers.  We refer the reader to \cite{HLL, HLLS, LWX, XLLM}.

Recently,  the averaging principle and LDP for  McKean-Vlasov SDEs driven by  fractional Brownian motions have been investigated.     Fan et al. \cite{FYY} studied the asymptotic behaviors for distribution dependent SDEs driven by fractional Brownian motions.  Shen et al. \cite{SXW} built the  the averaging principle for  distribution dependent SDEs driven by fractional Brownian motions. Aidara et al.  \cite{ASF} analyzed the averaging principle for BSDEs driven by two mutually independent fractional Brownian motions.  As far as we are concerned, there is no literature to discuss the LDP for slow-fast McKean-Vlasov SDEs driven by fractional  Brownian motions and Brownian motions.

In the paper, we shall consider the following slow-fast McKean-Vlasov SDEs driven by fractional  Brownian motions and Brownian motions:
\begin{align}\label{1}
\begin{cases}
 &\dif X^{\epsilon}_{t}=f_{1}(X^{\epsilon}_{t}, \sL_{X^{\epsilon}_{t}}, Y^{\epsilon}_{t} )\dif t + g_{1}(X^{\epsilon}_{t}, \sL_{X^{\epsilon}_{t} })\dif W^{1}_{t}+ l( \sL_{X^{\epsilon}_{t}})\dif B^{H}_{t} , X^{\epsilon}_{0}=x, \\
 & \dif Y^{\epsilon}_{t}=\frac{1}{\varepsilon}b(X^{\epsilon}_{t}, \sL_{X^{\epsilon}_{t}}, Y^{\epsilon}_{t} )\dif t + \frac{1}{\sqrt{\varepsilon}}\sigma_{1}(X^{\epsilon}_{t}, \sL_{X^{\epsilon}_{t}}, Y^{\epsilon}_{t} )\dif W^{1}_{t}\\
 &\quad \quad \quad + \frac{1}{\sqrt{\varepsilon}}\sigma_{2}(X^{\epsilon}_{t}, \sL_{X^{\epsilon}_{t}}, Y^{\epsilon}_{t} )\dif W^{2}_{t}, Y^{\epsilon}_{0}=y,
 \end{cases}
 \end{align}
where $b, f_{1}:\mR^{d}\times \sP_{2}(\mR^{d}) \times \mR^{d} \rightarrow \mR^{d}, $ $g_{1}: \mR^{d}\times \sP_{2}(\mR^{d})  \rightarrow \mR^{d\times d_{1}},$ $l:  \sP_{2}(\mR^{d})\rightarrow  \mR^{d\times d_{1}},$ $\sigma_{1}:\mR^{d}\times \sP_{2}(\mR^{d}) \times \mR^{d} \rightarrow \mR^{d\times d_{1}},$  $\sigma_{2}:\mR^{d}\times \sP_{2}(\mR^{d}) \times \mR^{d} \rightarrow \mR^{d\times d_{2}},$  $W^{1}, W^{2}$ are $d_{1}-$dimensional Brownian motion and $d_{2}-$Brownian motion, respectively, $B^{H}$ is a $d_{1}-$dimensional  fractional Brownian motion.

The main contributions are as follows:
\begin{itemize}
\item[$\bullet$] We give the definition of LDP  on the product space related to fractional  Brownian motion and Brownian motion, which is different from the traditional definition for LDP.

\item[$\bullet$] Two LDP criteria  are given on the product space related to fractional  Brownian motion and Brownian motion.

\item[$\bullet$]  The LDP is derived for  slow-fast McKean-Vlasov SDEs \eqref{1}  by using weak convergence method.
\end{itemize}

We close this part by giving our organization  for this article. In Section 2, we introduce some necessary notations, assumptions. In Section 3,  We give our main results. Throughout this paper, we make the following convention: the letter $C(\eta)$ with or without indices will denote different positive constants which only depends on $\eta$, whose value may
vary from one place to another.

\section{Preliminaries}
\subsection{Notations}Throughout this paper,  denote $C_0([0,T],\mR^d)$  by the continuous functions vanishing at $0$ equipped with the supremum norm. Let  $\Omega_{1}=C_{0}(0,T; \mR^{d_{1}}),\Omega_{2}=C_{0}(0,T; \mR^{d_{2}}), \Omega_{3}=C_{0}(0,T; \mR^{d_{3}}). $ $\sF_{i}, i=1,2,3$ are the Borel $\sigma-$algrebra, $\mF_{i}:=\{\sF^{i}_{t}, t\in [0,T]\}$ are the $\sigma-$algebra filtration  and $P^{i}, i=1,2,3$ are the probability on $\Omega_{i},, i=1,2,3$ such that
 $W^{1}_{\cdot}$ is a   $d_{1}$-Brownian motion on  $(\Omega_{1}, \sF_{1},\mF_{1}, P^{1}),$ $W^{2}_{\cdot}$ is a   $d_{2}$-Brownian motion on  $(\Omega_{2}, \sF_{2},\mF_{2}, P^{2})$ and $B^{H}_{\cdot}$ is a $d_{1}$-fractional Brownian motion $(\Omega_{3}, \sF_{3},\mF_{3}, P^{3}).$
  Let $(\Omega:=\Omega_{1}\times \Omega_{2} \times \Omega_{3},\sF:=\sF_{1}\times \sF_{2} \times \sF_{3}, \mF:=\mF_{1}\times \mF_{2} \times \mF_{3}, P:=P^{1}\times P^{2} \times P^{3}), $  where $\mF_{1}\times \mF_{2} \times \mF_{3}:= \{ \sF^{1}_{t}\times\sF^{2}_{t}\times \sF^{3}_{t} , t\in [0,T] \},$  be  the product space.
 Then, $W^{1}, W^{2},B^{H}$  are mutually independent.  If $x, y \in \mR^{d},$ we use $|x |$  to denote the Euclidean norm of
$x,$ and use  $\langle x, y\rangle$ or $xy$ to denote the Euclidean inner product. If  $A$ is a matrix, $A^{T}$  is the transpose  of $A,$  and   $|A |$ represents  $\sqrt{\mathrm{Tr} (AA^{T}).}$   Moreover, let $\lfloor a \rfloor$ be the integer parts of $a.$       Let $\sB(\mR^{d})$ be  the  Borel $\sigma-$algebra on $\mR^{d}$, $C(\mR^{d})$ denotes all continuous functions on $ \mR^{d}$ and
$C^{k}(\mR^{d})$ denotes all continuous functions on $ \mR^{d}$ with  continuous partial derivations
         of order up to $k$.
Let $\sP(\mR^{d})$ be the space  of all probability measures on $\sB(\mR^{d}) $, and  $\sP_{p}(\mR^{d})$ denotes the space  of all probability measures defined on $\sB(\mR^{d}) $ with finite $p$th moment:
$$[\mu(|\cdot|^{p})]^{\frac{1}{p}}:=\bigg(\bigg.\int_{\mR^{d}}|x|^{p}\mu(\dif x)\bigg)\bigg.^{\frac{1}{p}}<\infty.$$
For $ \mu, \nu\in \sP_{p}(\mR^{d})$, we define the
Wasserstein distance for $p\geq 1$  as follows: $$W_{p}(\mu, \nu):=\inf_{\pi \in \Pi(\mu, \nu)}\bigg\{\bigg. \int_{\mR^{d}\times \mR^{d}}|x-y|^{p}\pi(\dif x, \dif y)   \bigg\}\bigg.^{\frac{1}{p}}, $$
where $\Pi(\mu, \nu)$ is the family of all coupling for $\mu, \nu.$

\subsection{Fractional integral and derivative}
This aim of this section  is to introduce  some  notion  and notation of fractional calculus involving fractional integral and derivative,
Wiener space associated to fractional Brownian motion.

We mention  that $B^H=(B^{H,1},\cdots,B^{H,d})$ with Hurst parameter $H\in(0,1)$ is a centered Gaussian process with the covariance function $\mE(B^{H,i}_t B^{H,j}_s)=R_H(t,s)\delta_{i,j}$, where
\begin{align*}
 R_H(t,s)=\frac{1}{2}\left(t^{2H}+s^{2H}-|t-s|^{2H}\right),\ \  t,s\in[0,T].
\end{align*}
Then, the following results hold:
\begin{enumerate}
\item[$1^{\circ}$]
$\mE(|B_t^{H,i}-B_s^{H,i}|^q)=C_q|t-s|^{qH}$ for every $q\geq 1$ and $i=1,\cdots,d$.
\item[$2^{\circ}$]
 $B^H$ is $(H-\varepsilon)$-order H\"{o}lder continuous a.s. for any $\varepsilon\in(0,H)$ and is an $H$-self similar process.
\end{enumerate}
For any  $a,b\in\mR$ with $a<b,$
 $f\in L^1([a,b],\mR)$ and $\alpha>0$, the right-sided (respectively left-sided) fractional Riemann-Liouville integral of $f$ of order $\alpha$
on $[a,b]$ is defined as
\begin{align}\label{FrIn}
&I_{b-}^\alpha f(x)=\frac{(-1)^{-\alpha}}{\Gamma(\alpha)}\int_x^b\frac{f(y)}{(y-x)^{1-\alpha}}\dif y\\
&\qquad\left(\mbox{respectively}\ \ I_{a+}^\alpha f(x)=\frac{1}{\Gamma(\alpha)}\int_a^x\frac{f(y)}{(x-y)^{1-\alpha}}\dif y\right),\nonumber
\end{align}
where $x\in(a,b)$ a.e., $(-1)^{-\alpha}=e^{-i\alpha\pi}$ and $\Gamma$ denotes the Gamma function.

Fractional differentiation may be given as an inverse operation.
Let $\alpha\in(0,1)$ and $p\geq1$.
If $f\in I_{a+}^{\alpha}(L^p([a,b],\mR))$ (respectively $I_{b-}^\alpha(L^p([a,b],\mR)))$,   there exists a unique  function $g$ in $L^p([a,b],\mR)$ satisfying $f=I_{a+}^\alpha g$ (respectively $f=I_{b-}^\alpha g$)  and it coincides with the left-sided (respectively right-sided) Riemann-Liouville derivative
of $f$ of order $\alpha$ shown by
\begin{align*}
&D_{a+}^\alpha f(x)=\frac{1}{\Gamma(1-\alpha)}\frac{\dif}{\dif x}\int_a^x\frac{f(y)}{(x-y)^\alpha}\dif y\\
&\qquad\left(\mbox{respectively}\ D_{b-}^\alpha f(x)=\frac{(-1)^{1+\alpha}}{\Gamma(1-\alpha)}\frac{\dif}{\dif x}\int_x^b\frac{f(y)}{(y-x)^\alpha}\dif y\right).
\end{align*}
The corresponding Weyl representation is as follows:
\begin{align}\label{FrDe}
&D_{a+}^\alpha f(x)=\frac{1}{\Gamma(1-\alpha)}\left(\frac{f(x)}{(x-a)^\alpha}+\alpha\int_a^x\frac{f(x)-f(y)}{(x-y)^{\alpha+1}}\dif y\right)\\
&\qquad\left(\mbox{respectively}\ \ D_{b-}^\alpha f(x)=\frac{(-1)^\alpha}{\Gamma(1-\alpha)}\left(\frac{f(x)}{(b-x)^\alpha}+\alpha\int_x^b\frac{f(x)-f(y)}{(y-x)^{\alpha+1}}\dif y\right)\right).\nonumber
\end{align}
Obviously, from \cite{SK},  the convergence of the integrals at the singularity $y=x$ holds pointwise for almost all $x$ if $p=1$ and in the $L^p$ sense if $p>1$.

\subsection{Wiener space associated to fractional Brownian motion}

Denote by  $\mathscr{E}$  the set of step functions on $[0,T]$ and the Hilbert space $\mathcal {H}$ is  the closure of
$\mathscr{E}$ with respect to the scalar product
\begin{align*}
\left\langle (\mathbb I_{[0,t_1]},\cdot\cdot\cdot,\mathbb I_{[0,t_{d_{1}}]}),(\mathbb I_{[0,s_1]},\cdot\cdot\cdot,\mathbb I_{[0,s_{d_{1}}]})\right\rangle_{\mathcal{H}}=\sum_{i=1}^{d_{1}} R_H(t_i,s_i).
\end{align*}
In order to give the  integral representation of $R_H(t,s), $   we set
\begin{align*}
K_H(t,s):=\Gamma\left(H+\frac{1}{2}\right)^{-1}(t-s)^{H-\frac{1}{2}}F\left(H-\frac{1}{2},\frac{1}{2}-H,H+\frac{1}{2},1-\frac{t}{s}\right),
\end{align*}
in which $F(\cdot,\cdot,\cdot,\cdot)$ is the Gauss hypergeometric function ( see \cite{Decreusefond&Ustunel98a} or \cite{Nikiforov&Uvarov88}). Obviously, $K_H$ is a square integrable kernel.
Then $R_H(t,s)$ has the following integral representation from \cite{Decreusefond&Ustunel98a}:
\begin{align*}
 R_H(t,s)=\int_0^{t\wedge s}K_H(t,r)K_H(s,r)\dif r.
\end{align*}
Furthermore, the mapping $(\mathbb I_{[0,t_1]},\cdot\cdot\cdot,\mathbb I_{[0,t_{d_{1}}]})\mapsto\sum_{i=1}^dB_{t_i}^{H,i}$ can be extended to an isometry between $\mathcal{H}$  and the Gaussian space $\mathcal {H}_1$ associated with $B^H$.
Denote this isometry by $\varphi\mapsto B^H(\varphi)$. Consequently, due to \cite{Alos&Mazet&Nualart01a} again,
 $B^H$ has the following Volterra-type representation
\begin{align}\label{IRFor}
B_t^H=\int_0^tK_H(t,s)\dif W_s, \ \ t\in[0,T],
\end{align}
where $W$  is some $d_{1}$-dimensional Wiener process defined on $(\Omega_{3},\sF_{3},P^{3}).$

Besides,  define the operator $K_H: L^2([0,T],\mathbb{R}^{d_{1}})\rightarrow I_{0+}^{H+1/2 }(L^2([0,T],\mathbb{R}^{d_{1}}))$ by
\begin{align*}
 (K_H f)(t)=\int_0^tK_H(t,s)f(s)\dif s.
\end{align*}
According  to \cite{Decreusefond&Ustunel98a}, we obtain   that it is an isomorphism and for each $f\in L^2([0,T],\mathbb{R}^{d_{1}})$,
\begin{equation}\nonumber
(K_H f)(s)=
\left\{
\begin{array}{ll}\vspace{0.3cm}
I_{0+}^{2H}s^{1/2-H}I_{0+}^{1/2-H}s^{H-1/2}f,\ \ H\in(0,1/2),\\
I_{0+}^{1}s^{H-1/2}I_{0+}^{H-1/2}s^{1/2-H}f,\ \ H\in(1/2,1).
\end{array} \right.
\end{equation}
Thus,  for any  $g\in I_{0+}^{H+1/2}(L^2([0,T],\mR^{d_{1}}))$, the inverse operator of $K_H$ is given by
\begin{equation}\label{InOp}
(K_H^{-1}g)(s)=
\left\{
\begin{array}{ll}\vspace{0.3cm}
s^{1/2-H}D_{0+}^{1/2-H}s^{H-1/2}D_{0+}^{2H}g,\ \ H\in(0,1/2),\\
s^{H-1/2}D_{0+}^{H-1/2}s^{1/2-H}g',\ \ \ \ \ \ \ H\in(1/2,1).
\end{array} \right.
\end{equation}
In particular, if $g$ is absolutely continuous, we have
\begin{align}\label{-HOF}
(K_H^{-1}g)(s)=s^{H-\frac{1}{2}}I_{0+}^{\frac{1}{2}-H}s^{\frac{1}{2}-H}g',\ \ H\in(0,1/2).
\end{align}

Now, define the linear operator $K_H^*:\mathscr{E}\rightarrow L^2([0,T],\mR^{d_{1}})$ as follows
\begin{align*}
(K_H^*\varphi)(s)=K_H(T,s)\varphi(s)+\int_s^T(\varphi(r)-\varphi(s))\frac{\partial K_H}{\partial r}(r,s)\dif r.
\end{align*}
According to \cite{Alos&Mazet&Nualart01a},
 we have the following the relation  for any $\varphi,\phi\in\mathscr{E},$
\begin{align}\label{-+}
\langle K_H^*\varphi,K_H^*\phi\rangle_{L^2([0,T],\mR^{d_{1}})}=\langle\varphi,\phi\rangle_\mathcal {H}=H(2H-1)\int^{T}_{0}\int^{T}_{0}|t-s|^{2H-1}\langle \varphi(s), \phi(t)\rangle\dif s\dif t.
\end{align}
and then  the bounded linear transform theorem implies that  $K_H^*$ can be extended to an isometry between $\mathcal{H}$ and $L^2([0,T],\mR^{d_{1}})$.

Note that the injection $R_H=K_H\circ K_H^*:\mathcal{H}\rightarrow\Omega_{3}$ embeds $\mathcal{H}$ densely into $\Omega_{3}$ and
for every $\psi\in\Omega_{3}^*\subset\mathcal{H}$ it  holds that
$\mE e^{i\langle B^H,\psi\rangle}=\exp(- 1 2\|\psi\|^2_\mathcal{H})$.
By (2.6) in \cite{FYY}, for any $\varphi \in \mathcal{H},$ we derive that
\begin{align}\label{+-}
(R_{H}\varphi)(t)= \int^{t}_{0}\int^{s}_{0}\frac{\partial K_{H}}{\partial s}(s,r)(K^{*}_{H}\varphi)(r)\dif r \dif s.
\end{align}

Next, we introduce some results about  the Malliavin calculus for fractional Brownian motion.
Let  $\sS$ be the set of smooth and cylindrical random variables of the form
\begin{align*}
F=f(B^H(\varphi_1),\cdot\cdot\cdot,B^H(\varphi_n)),
\end{align*}
where $n\geq 1, f\in C_b^\infty(\mathbb{R}^n)$ and
$\varphi_i\in\mathcal{H}, 1\leq i\leq n$.
The Malliavin derivative of $F$  is given as follows:
\begin{align*}
\mD F=\sum_{i=1}^n\frac{\partial f}{\partial x_i}(B^H(\varphi_1),\cdot\cdot\cdot,B^H(\varphi_n))\varphi_i,
\end{align*}
and denoted by $\mD F.$

For any $p\geq 1$, we denote by $\mathbb{D}^{1,p}$  the Sobolev space which is  the completion of $\mathcal {S}$ with respect to the norm
\begin{align*}
\|F\|_{1,p}^p=\mathbb{E}|F|^p+\mathbb{E}\|\mD F\|^p_{\mathcal {H}}.
\end{align*}
Denote by $\delta$ and $\mathrm{Dom}\delta$ the dual operator of $\mathbb{D}$ and its domain, respectively.
The following results are needed later. By \cite[Proposition 5.2.1]{ND} and \cite[Proposition 5.2.2]{ND}, one has
\begin{prp}\label{Tra-Prp1}
Denote by $\mD^W$  the derivative operator with respect to the underlying Wiener process $W$  in \eqref{IRFor}, and $\mD_W^{1,2}$ the corresponding Sobolev space. For every  $F\in\mD_W^{1,2}=\mD^{1,2}$, we obtain that
\begin{align*}
K_H^*\mD F=\mD^WF,
\end{align*}
 \end{prp}

\begin{prp}\label{Tra-Prp2}
$\mathrm{Dom}\delta=(K_H^*)^{-1}(\mathrm{Dom}\delta_W)$,
and it holds that  $\delta(u)=\delta_W(K_H^*u)$ for any $\mathcal{H}$-valued random variable $u$ in $\mathrm{Dom}\delta$,
where $\delta_W$ represents the divergence operator corresponding  to the underlying Wiener process $W$  in \eqref{IRFor}.
\end{prp}

\begin{rem}\label{Tra-Rem1}
Combining the above proposition and  \cite[Proposition 1.3.11]{ND}, it implies  that $u\in\mathrm{Dom}\delta$
if $K_H^*u\in L_a^2([0,T]\times\Omega,\mR^{d_{1}})$ (the closed subspace of $L^2([0,T]\times\Omega,\mR^{d_{1}})$ formed by the adapted processes).
\end{rem}

Finally,  we complete this section by giving the following  notation for future use.
\begin{align*}
&\mathcal{A}:=\bigg\{\bigg. \bar{h}\,\mbox{is }\, \mR^{d_{1}}\, \mbox{valued}\, \sF_{t}-\mbox{predictable process such that }  \|\bar{h}\|_{\mathcal{H}}< \infty\bigg\}\bigg..\\
&\mathcal{S}_{N}:=\{\bar{h}\in \mathcal{H}: \|\bar{h}\|_{\mathcal{H}}^{2} \leq N\}.\\
&\mathcal{A}_{N}:=\bigg\{\bigg.\bar{h}\in \mathcal{A}: \bar{h}(\omega)\in  \mathcal{S}_{N} \bigg\}\bigg..\\
&\mH:=\bigg\{\bigg. h=\int^{\cdot}_{0}\dot{h}(s)\dif s: \dot{h} \in L^{2}(0,T;\mR^{d_{1}+d_{2}})\,\mbox{with the norm}\\
&~~~~~~~~~~~~~~~~~~~~~~~~~~~~~~~ \|h\|_{\mH}:=\bigg(\bigg. \int^{T}_{0}|\dot{h}(s)|^{2}\dif s \bigg)\bigg.^{\frac{1}{2}} < \infty \bigg\}\bigg..\\
&\mA:=\bigg\{\bigg. h\,\mbox{is }\, \mR^{d_{1}+d_{2}}\, \mbox{valued}\, \sF_{t}-\mbox{predictable process such that }  \int^{T}_{0}|\dot{h}(s)|^{2}\dif s < \infty\bigg\}\bigg..\\
&\mS_{N}:=\{h\in \mH:  \int^{T}_{0}|\dot{h}(s)|^{2}\dif s \leq N\}.\\
&\mA_{N}:=\bigg\{\bigg.h\in \mA: h(\omega)\in  \mS_{N} \bigg\}\bigg..
 \end{align*}
For $h_{\epsilon}\in \mA_{N}, h\in \mA_{N}, $ $h_{\epsilon}\Rightarrow h$ means $h_{\epsilon}$ converges in distribution to   $h$  as $\epsilon \rightarrow 0.$\\
For $\bar{h}_{\epsilon}\in \mathcal{A}_{N}, \bar{h}\in \mathcal{A}_{N}, $ $\bar{h}_{\epsilon}\Rightarrow \bar{h}$ means $\bar{h}_{\epsilon}$ converges in distribution to   $\bar{h}$  as $\epsilon \rightarrow 0.$

\subsection{Assumptions}

We now impose the following condition on coefficients:
\begin{itemize}
\item[(H1)] Let $C$ and $\alpha$ be two positive constants such that the following conditions hold for any $ x, x_{1}, x_{2}\in \mR^{d}, \mu, \mu_{1}, \mu_{2}\in \sP_{2}(\mR^{d}),  y, y_{1}, y_{2}\in \mR^{d}:$
\begin{align}\label{01}
|&f_{1}(x_{1},\mu_{1}, y_{1})-f_{1}(x_{2},\mu_{2}, y_{2})|
+|b(x_{1},\mu_{1}, y_{1})-b(x_{2},\mu_{2}, y_{2})| \no\\
 &+|\sigma_{1}(x_{1},\mu_{1}, y_{1})-\sigma_{1}(x_{2},\mu_{2}, y_{2})|
+ |\sigma_{2}(x_{1},\mu_{1}, y_{1})-\sigma_{2}(x_{2},\mu_{2}, y_{2})|\no\\
& \leq C(|t_{1}-t_{2}|+|x_{1}-x_{2}|+|y_{1}-y_{2}|+W_{2}(\mu_{1}, \mu_{2})),\\
 &|g_{1}(x_{1},\mu_{1})-g_{1}(x_{2},\mu_{2})|
 \leq C(|x_{1}-x_{2}|+W_{2}(\mu_{1}, \mu_{2})), \label{01a}\\
& |l(\mu_{1})-l(\mu_{2})||
 \leq CW_{2}(\mu_{1}, \mu_{2}), \label{01b}
\end{align}
and
\begin{align}\label{03}
4\langle &b(x,\mu, y_{1})-b(x,\mu, y_{2}),  y_{1}-y_{2}  \rangle + 6|\sigma_{1}(x,\mu, y_{1})-\sigma_{1}(x,\mu, y_{2})|^{2}\no\\
&+ 6|\sigma_{2}(x,\mu, y_{1})-\sigma_{2}(x,\mu, y_{2})|^{2}\leq -\alpha | y_{1}-y_{2}|^{2}.
\end{align}

\end{itemize}

Firstly, we give several uniform estimations $w.r.t.$ $\varepsilon \in (0, 1) $ for the 4th moment of solution $(X^{\epsilon}, Y^{\epsilon})$ of Eq.  \eqref{1}.
\bl \label{1+}
Assume $(\mathrm{H1})$. Then the following inequalities hold:
\begin{align*}
\sup_{\varepsilon \in (0,1)}\sup_{t\in [0, T]}\mE[|X^{\epsilon}_{t}|^{4}]\leq C(T)(1+|x|^{4}+|y|^{4}),
\end{align*}
\begin{align*}
\sup_{\varepsilon \in (0,1)}\sup_{t\in [0, T]}\mE[|Y^{\epsilon}_{t}|^{4}]\leq C(T)(1+|x|^{4}+|y|^{4}).
\end{align*}
\el

\begin{proof}
By Eq.  \eqref{1}, we have
\begin{align}\label{2}
\mE[\sup_{0\leq s \leq t}|X_s^\epsilon|^{4}]&\leq C_{1}|x|^{4} +C_{1}\mE\bigg[\bigg. \sup_{0\leq s \leq t}\bigg|\bigg. \int^{s}_{0}f_{1}(X^{\epsilon}_{r}, \sL_{X^{\epsilon}_{r}}, Y^{\epsilon}_{r} )\dif r \bigg|\bigg.^{4} \bigg]\bigg. \no\\
& \quad+C_{1}\mE\bigg[\bigg. \sup_{0\leq s \leq t}\bigg|\bigg. \int^{s}_{0}g_{1}(X^{\epsilon}_{r}, \sL_{X^{\epsilon}_{r}}, Y^{\epsilon}_{r}  ) \dif W_{r} \bigg|\bigg.^{4} \bigg]\bigg.\no\\
&\quad + C_{1}\mE\bigg[\bigg. \sup_{0\leq s \leq t}\bigg|\bigg. \int^{s}_{0}l(r, \sL_{X^{\epsilon}_{r}})  \dif B^{H}_{r} \bigg|\bigg.^{4} \bigg]\bigg.\no\\
&=: C_{1}|x|^{4}+ I_{1}+I_{2}+I_{3}.
\end{align}
where $ C_{1}$ is a constant.

For the term $I_{1},$   by the fact $W_{2}^{4}(\sL_{X_{t}}, \delta_{0})\leq \mE[|X_{t}|^{4}]$    and $(\mathrm{H1})$, we derive
\begin{align*}
 I_{1}&=C_{1}\mE\bigg[\bigg. \sup_{0\leq s \leq t}\bigg|\bigg. \int^{s}_{0}f_{1}(X^{\epsilon}_{r}, \sL_{X^{\epsilon}_{r}}, Y^{\epsilon}_{r} )] \dif r \bigg|\bigg.^{4} \bigg]\bigg. \\
 & \leq C(T)\mE\int^{t}_{0}(|X^{\epsilon}_{s}|^{4}+|Y^{\epsilon}_{s}|^{4})\dif s+C(T).
\end{align*}
Next, we look at $I_{2}.$ The Burkholder-Davis-Gundy(BDG) inequality and $(\mathrm{H1})$ yield that
\begin{align*}
 I_{2}&=C_{1}\mE\bigg[\bigg. \sup_{0\leq s \leq t}\bigg|\bigg. \int^{s}_{0}g_{1}(X^{\epsilon}_{r}, \sL_{X^{\epsilon}_{r}} )  \dif W_{r} \bigg|\bigg.^{4} \bigg]\bigg. \\
 & \leq C(T)\mE\bigg[\bigg.\bigg|\bigg.\int^{t}_{0}|g_{1}(X^{\epsilon}_{t}, \sL_{X^{\epsilon}_{r}} )|^{2}\dif s\bigg|\bigg.^{2}\bigg]\bigg.\\
& \leq C(T)\mE\int^{t}_{0}|X^{\epsilon}_{s}|^{4}\dif s+ C(T).
\end{align*}
Finally,  for $I_{3},$ it follows from \cite [Theorem 3.2]{AN},  (3.5) in \cite{FYY} and  $(\mathrm{H1})$,  we obtain
\begin{align*}
 I_{3}&=C_{1}\mE\bigg[\bigg. \sup_{0\leq s \leq t}\bigg|\bigg. \int^{s}_{0}l(s, \sL_{X^{\epsilon}_{r}})  \dif B^{H}_{r} \bigg|\bigg.^{4} \bigg]\bigg. \\
 & \leq C(T, H)\mE\bigg[\bigg.\bigg|\bigg.\int^{t}_{0}|l(s, \sL_{X^{\epsilon}_{s}})|^{2}\dif s\bigg|\bigg.^{2}\bigg]\bigg.\\
& \leq C(T, H)\mE\int^{t}_{0}|X^{\epsilon}_{s}|^{4}\dif s+C(T).
\end{align*}
These imply
\begin{align*}%\label{3}
\mE[\sup_{0\leq s \leq t}|X^{\epsilon}_s|^{4}]\leq C_{1}|x|^{4}+C(T) +C(T, H)\int^{t}_{0}(\mE[\sup_{0\leq r \leq s}|X^{\epsilon}_{r}|^{4}]+\mE[\sup_{0\leq r \leq s}|Y^{\epsilon}_{r}|^{4}])\dif s.
\end{align*}
An application of  It\^{o}'s formula yields that
\begin{align}\label{3}
\mE[|Y^{\epsilon}_{t}|^{4}]&=|y|^{4}+\frac{4}{\varepsilon}\mE\int^{t}_{0}|Y^{\epsilon}_{s}|^{2}\langle Y^{\epsilon}_{s},   b(s,X^{\epsilon}_{s}, \sL_{X^{\epsilon}_{s}}, Y^{\epsilon}_{s} )\rangle \dif s\no\\
& \quad+ \frac{2}{\varepsilon}\mE\int^{t}_{0}|Y^{\epsilon}_{s}|^{2}| \sigma_{1}(s,X^{\epsilon}_{s}, \sL_{X^{\epsilon}_{s}}, Y^{\epsilon}_{s} )|^{2} \dif s\no\\
& \quad+ \frac{4}{\varepsilon}\mE\int^{t}_{0}| \langle Y^{\epsilon}_{s},   \sigma_{1}(s,X^{\epsilon}_{s}, \sL_{X^{\epsilon}_{s}}, Y^{\epsilon}_{s} )\rangle |^{2} \dif s\no\\
& \quad+ \frac{2}{\varepsilon}\mE\int^{t}_{0}|Y^{\epsilon}_{s}|^{2}| \sigma_{2}(s,X^{\epsilon}_{s}, \sL_{X^{\epsilon}_{s}}, Y^{\epsilon}_{s} )|^{2} \dif s\no\\
& \quad+ \frac{4}{\varepsilon}\mE\int^{t}_{0}| \langle Y^{\epsilon}_{s},   \sigma_{2}(s,X^{\epsilon}_{s}, \sL_{X^{\epsilon}_{s}}, Y^{\epsilon}_{s} )\rangle |^{2} \dif s.
\end{align}
By \eqref{03}, there exists $\alpha'>0$ such that for any $t\in [0, T],$
\begin{align*}
\frac{\dif }{\dif t}\mE[|Y^{\epsilon}_{t}|^{4}]&\leq \frac{1}{\varepsilon}\mE[4|Y^{\epsilon}_{t}|^{2}\langle Y^{\epsilon}_{t},   b(s,X^{\epsilon}_{t}, \sL_{X^{\epsilon}_{t}}, Y^{\epsilon}_{t} )\rangle  + 6|Y^{\epsilon}_{s}|^{2}| \sigma_{1}(s,X^{\epsilon}_{s}, \sL_{X^{\epsilon}_{s}}, Y^{\epsilon}_{s} )|^{2}]\no\\
& \quad + 6|Y^{\epsilon}_{s}|^{2}| \sigma_{2}(s,X^{\epsilon}_{s}, \sL_{X^{\epsilon}_{s}}, Y^{\epsilon}_{s} )|^{2}]\no\\
& \quad \leq-\frac{\alpha'}{\varepsilon}\mE[|Y^{\epsilon}_{t}|^{4}] +\frac{C(T)}{\varepsilon}(\mE[|X^{\epsilon}_{t}|^{4}] +1).
\end{align*}
It follows from the comparison theorem that
\begin{align*}
\mE[|Y^{\epsilon}_{t}|^{4}]&\leq |y|^{4}e^{-\frac{\alpha' t}{\varepsilon}}+\frac{C(T)}{\varepsilon}\int^{t}_{0}e^{-\frac{\alpha' (t-s)}{\varepsilon}}(\mE[|X^{\epsilon}_{s}|^{4}] +1)\dif s\no\\
& \quad \leq |y|^{4} +C(T)(\sup_{0\leq s\leq t}\mE[|X^{\epsilon}_{t}|^{4}] +1).
\end{align*}
This together with \eqref{3} implies
\begin{align*}
\sup_{0\leq s\leq t}\mE[|X^{\epsilon}_{t}|^{4}]
 \leq C(T)|y|^{4} +C(T)\int^{t}_{0}\sup_{0\leq r\leq s}\mE[|X^{\epsilon}_{r}|^{4}]\dif s.
\end{align*}
Gronwall's  inequality gives
\begin{align*}
\sup_{0<\varepsilon< 1}\sup_{0\leq s\leq t}\mE[|X^{\epsilon}_{t}|^{4}]
 \leq C(T)(1+|x|^{4}+|y|^{4}).
\end{align*}
which also implies
\begin{align*}
\sup_{0<\varepsilon< 1}\sup_{0\leq s\leq t}\mE[|Y^{\epsilon}_{t}|^{4}]
 \leq C(T)(1+|x|^{4}+|y|^{4}).
\end{align*}
\end{proof}

\bl\label{5+}
Assume $(\mathrm{H1}).$ For any $0\leq t \leq t+u \leq  T,$  we have the following inequality.
$$ \mE[|X^{\varepsilon}_{t+u}-X^{\varepsilon}_{t}|^{2}] \leq C(T)(1+|x|^{2}+|y|^{2})u.$$
\el
\begin{proof}
By \eqref{1},  we have
\begin{align*}
 X^{\epsilon}_{t+u}- X^{\epsilon}_{t}&=\int^{t+u}_{t}f_{1}(X^{\epsilon}_{s}, \sL_{X^{\epsilon}_{s}}, Y^{\epsilon}_{s} ) \dif s \no\\
 &\quad + \int^{t+u}_{t}g_{1}(X^{\epsilon}_{s}, \sL_{X^{\epsilon}_{s}} ) \dif W^{1}_{s}+\int^{t+u}_{t} l( \sL_{X^{\epsilon}_{s}})\dif B^{H}_{s}.
\end{align*}
Set $\hat{X}_{t}:= X^{\epsilon}_{t+u}- X^{\epsilon}_{t}.$  Analogous to the calculation of \eqref{2}, it holds that
\begin{align*}
&\mE[|\hat{X}(t)|^{2}]\leq C\mE\bigg[\bigg. \bigg|\bigg. \int^{t+u}_{t}f_{1}(X^{\epsilon}_{r}, \sL_{X^{\epsilon}_{r}}, Y^{\epsilon}_{r} )  \dif r \bigg|\bigg.^{2} \bigg]\bigg. \no\\
& \quad+C\mE\bigg[\bigg. \bigg|\bigg. \int^{t+u}_{t}g_{1}(X^{\epsilon}_{r}, \sL_{X^{\epsilon}_{r}}, Y^{\epsilon}_{r} )   \dif W_{r} \bigg|\bigg.^{2} \bigg]\bigg.
 + C\mE\bigg[\bigg. \bigg|\bigg. \int^{t+u}_{t}l( \sL_{X^{\epsilon}_{r}})  \dif B^{H}_{r} \bigg|\bigg.^{2} \bigg]\bigg.\no\\
& \leq  C(T)u\int^{t+u}_{t}(\mE[1+ |X^{\epsilon}_{s}|^{2}]+\mE[|Y^{\epsilon}_{s}|^{2}])\dif s +  C(T)\int^{t+u}_{t}(\mE[1+ |X^{\epsilon}_{s}|^{2}]+\mE[|Y^{\epsilon}_{s}|^{2}])\dif s.
\end{align*}
 Then, Lemma \ref{1+} leads to the required assertion.
\end{proof}

\subsection{The average equation.} We now introduce the following  parameterized McKean-Vlasov equation: for fixed $t\geq 0, x\in \mR^{d}, \mu\in \sP_{2}(\mR^{d}),$ let

\begin{align}\label{8}
\dif Y^{x,\mu,y}_{s}&=b(x, \mu, Y^{x,\mu,y}_{s} )\dif s + \sigma_{1}(x, \mu, Y^{x,\mu,y}_{s} )\dif \tilde{W}^{1}_{s}\no\\
 &\quad  + \sigma_{2}(x, \mu, Y^{x,\mu,y}_{s} )\dif \tilde{W}^{2}_{s}, ~~ Y^{x,\mu,y}_{0}= y,
\end{align}
where $ \widetilde{W}^{1}_{\cdot},  \widetilde{W}^{2}_{\cdot}$ are $d_{1}$-dimensional and $d_{2}$-dimensional Brownian motions, respectively, on another complete probability space  $(\widetilde{\Omega}, \widetilde{\sF}, \widetilde{P})$ and $\{\widetilde{\sF}_{t}\}_{t\geq 0}$ is the natural filtration  generated by $\widetilde{W}^{1}_{t}, \widetilde{W}^{2}_{t}, B^{H}_{t}.$  Under $(\mathrm{H1})$, Eq. \eqref{8} has a unique strong solution   $\{Y^{t,x,\mu,y}_{s} \}_{s\geq 0}$ and it is a homogeneous  Markov process with the following estimate
$$\sup_{s\geq 0}\widetilde{\mE }[|Y^{x,\mu,y}_{s}|^{2}]\leq C(T)(1+|x|^{2}+|y|^{2}+\mu(|\cdot|^{2})).$$

Let $\{P^{x, \mu}_{s}\}_{s\geq 0}$ be the transition semigroup of $Y^{x,\mu,y}_{s}.$  By \cite [Theorem 4.3.9]{LR}, under $(\mathrm{H1}),$    $\{P^{x, \mu}_{s}\}_{s\geq 0}$ has a unique invariant  measure $\nu^{x,\mu}$ satisfying
$$\int_{\mR^{m}} |y|\nu^{x,\mu}(\dif y )\leq C(T)(1+|x|+|y|+[\mu(|\cdot|^{2})]^{\frac{1}{2}}). $$
Set

$$\bar{f}(x, \mu):=\int_{\mR^{d}}f_{1}( x, \mu, z)\nu^{x,\mu}(\dif z ).$$

We have the following lemma.
\bl\label{14}
$\bar{f}$  satisfies  the following Lipschitz conditions, i.e.
$$ |\bar{f}(x_{1}, \mu_{1} )- \bar{f}(x_{2}, \mu_{2} )| \leq C(|x_{1}-x_{2} |+W_{2}(\mu_{1},\mu_{2} )).$$
\el
\begin{proof}

By the definition of $\bar{f}$ and the Lipschitz continuity of $f,$  we derive
\begin{align*}
&|\bar{f}(x_{1}, \mu_{1})-\bar{f}( x_{2}, \mu_{2})|\\
&=\bigg|\bigg.\int_{\mR^{d_{1}}}f_{1}( x_{1}, \mu_{1}, z)\nu^{ x_{1}, \mu_{1}}(\dif z)- \int_{\mR^{d_{1}}}f_{1}(t, x_{2}, \mu_{2}, z)\nu^{ x_{2}, \mu_{2}}(\dif z)   \bigg|\bigg.\\
&\leq \bigg|\bigg.\int_{\mR^{d_{1}}}[f_{1}( x_{1}, \mu_{1}, z)- f_{1}( x_{2}, \mu_{2}, z)]\nu^{x_{1}, \mu_{1}}(\dif z)   \bigg|\bigg.\\
& \quad +\bigg|\bigg.\int_{\mR^{d_{1}}}f_{1}( x_{2}, \mu_{2}, z)\nu^{ x_{1}, \mu_{1}}(\dif z)- \int_{\mR^{d_{1}}}f_{1}( x_{2}, \mu_{2}, z)\nu^{ x_{2}, \mu_{2}}(\dif z)   \bigg|\bigg.\\
&\leq C(|x_{1}-x_{2}|+W_{2}^{2}(\mu_{1},\mu_{2})+W_{2}^{2}(\nu^{ x_{1}, \mu_{1}},\nu^{ x_{2}, \mu_{2}})).
\end{align*}
By \cite[Theorem 3.1]{WFY},  we obtain
\begin{align*}
&W_{2}^{2}(\nu^{ x_{1}, \mu_{1}},\nu^{ x_{2}, \mu_{2}})\\
&\leq 3W_{2}^{2}(\nu^{ x_{1}, \mu_{1}}, \sL_{Y^{x,\mu_{1}, 0}_{s}})+3W_{2}^{2}(\nu^{ x_{2}, \mu_{2}}, \sL_{Y^{x,\mu_{2}, 0}_{s}})+ 3W_{2}^{2}(\sL_{Y^{x,\mu_{1}, 0}_{s}}, \sL_{Y^{x,\mu_{2}, 0}_{s}})\\
& \leq Ce^{-2\lambda_{0}s}(W_{2}^{2}(\nu^{ x_{1}, \mu_{1}}, \delta_{0})+W_{2}^{2}(\nu^{ x_{1}, \mu_{1}}, \delta_{0}) )+C(|x_{1}-x_{2}|+W_{2}^{2}(\mu_{1},\mu_{2})).
\end{align*}
From the above calculation,  we obtain
$$ |\bar{f}(x_{1}, \mu_{1} )- \bar{f}(x_{2}, \mu_{2} )|\leq C(|x_{1}-x_{2} |+W_{2}(\mu_{1},\mu_{2} )). $$
\end{proof}

 We now consider the following average equation:
\begin{align}\label{001}
\dif \bar{X}_{t}&=\bar{f}(\bar{X}_{t}, \sL_{\bar{X}_{t}} )\dif t +g_{1}(\bar{X}_{t}, \sL_{\bar{X}_{t}} )\dif W^{1}_{t}+ l( \sL_{\bar{X}_{t}})\dif B^{H}_{t} , \bar{X}_{0}=x.
\end{align}

Since $\bar{f}_{1}, \bar{g}_{1}, l$ satisfy the Lipschitz condition, we have the following existence and uniqueness result for  Eq. \eqref{001}.
\bt

Assume  $(\mathrm{H}1).$  Then Eq. \eqref{001} has a unique solution.

\et

\subsection{Large deviation principle }
Consider the following multi-scale McKean-Vlasov system with small perturbation.
\begin{align}\label{31}
\begin{cases}
 &\dif X^{\epsilon,\varepsilon}_{t}=f_{1}(X^{\epsilon,\varepsilon}_{t}, \sL_{X^{\epsilon,\varepsilon}_{t}}, Y^{\epsilon,\varepsilon}_{t} )\dif t + \epsilon^{H}l( \sL_{X^{\epsilon,\varepsilon}_{t}})\dif B^{H}_{t} \\
 &\quad \quad \quad+\sqrt{\epsilon}g_{1}(X^{\epsilon,\varepsilon}_{t}, \sL_{X^{\epsilon,\varepsilon}_{t} }  )\dif W^{1}_{t},  X^{\epsilon,\varepsilon}_{0}=x_{0}, \\
 & \dif Y^{\epsilon,\varepsilon}_{t}=\frac{1}{\varepsilon}b(X^{\epsilon,\varepsilon}_{t}, \sL_{X^{\epsilon,\varepsilon}_{t}}, Y^{\epsilon,\varepsilon}_{t} )\dif t + \frac{1}{\sqrt{\varepsilon}}\sigma_{1}(X^{\epsilon,\varepsilon}_{t}, \sL_{X^{\epsilon,\varepsilon}_{t}}, Y^{\epsilon,\varepsilon}_{t} )\dif W^{1}_{t}\\
 &\quad \quad \quad + \frac{1}{\sqrt{\varepsilon}}\sigma_{2}(X^{\epsilon,\varepsilon}_{t}, \sL_{X^{\epsilon,\varepsilon}_{t}}, Y^{\epsilon,\varepsilon}_{t} )\dif W^{2}_{t}, Y^{\epsilon,\varepsilon}_{0}=y_{0}.
 \end{cases}
 \end{align}
Thus, due to the classical Yamada-Watanabe theorem,  there exists a measurable map $\Gamma_{\sL_{X^{\epsilon,\varepsilon}_{t}}}: C([0,T];\mR^{d_{1}+d_{2}})\times C(0,T; \mR^{d_{1}})\rightarrow C(0,T; \mR^{d})$ such that we have the following representation
$$X^{\epsilon,\varepsilon}_{t}=\Gamma_{\sL_{X^{\epsilon,\varepsilon}_{t}}}(\sqrt{\epsilon}W, \epsilon^{H}B^{H}),$$
where $W:=(W^{1},W^{2}).$
For simplicity of notation,  we denote $\Gamma^{\epsilon}=\Gamma_{\sL_{X^{\epsilon,\varepsilon}_{t}}}.$
For $\tilde{h}^{\epsilon}=(h^{\epsilon}, \bar{h}^{\epsilon} ),  h^{\epsilon}\in \mA_{N}, \bar{h}^{\epsilon}\in \mathcal{A}_{N},$    $X^{\epsilon, \varepsilon, \tilde{h}^{\epsilon}}:=\Gamma^{\epsilon}\bigg(\bigg.\sqrt{\epsilon}W_{\cdot}+\int^{\cdot}_{0} \dot{h}^{\epsilon}(\cdot)\dif s, \epsilon^{H}B^{H}_{\cdot}+ R_{H}\bar{h}^{\epsilon}(\cdot)   \bigg)\bigg.$ is the first part of solution of the following equation.
\begin{align}\label{42}
\begin{cases}
&\dif X^{\epsilon, \varepsilon, \tilde{h}^{\epsilon}}_{t}=f_{1}(X^{\epsilon, \varepsilon, \tilde{h}^{\epsilon}}_{t}, \sL_{X^{\epsilon, \varepsilon}_{t}}, Y^{\epsilon, \varepsilon, \tilde{h}^{\epsilon}}_{t} )\dif t+  g_{1}(X^{\epsilon, \varepsilon, \tilde{h}^{\epsilon}}_{t}, \sL_{X^{\epsilon, \varepsilon}_{t} }  ) {\cal P}_{1}\dot{h}^{\epsilon}(t)\dif t
  +l( \sL_{X^{\epsilon, \varepsilon}_{t}})\dif (R_{H}\bar{h}^{\epsilon})(t)\\
 &\quad \quad \quad+\sqrt{\epsilon}g_{1}(X^{\epsilon, \varepsilon, \tilde{h}^{\epsilon}}_{t}, \sL_{X^{\epsilon, \varepsilon}_{t} }  ) \dif W^{1}_{t}
   + \epsilon^{H}l( \sL_{X^{\epsilon, \varepsilon}_{t}})\dif B^{H}_{t} , X^{\epsilon}_{0}=x, \\
 &  \dif Y^{\epsilon, \varepsilon, \tilde{h}^{\epsilon}}_{t}=\frac{1}{\varepsilon}b(X^{\epsilon, \varepsilon, \tilde{h}^{\epsilon}}_{t}, \sL_{X^{\epsilon, \varepsilon}_{t}}, Y^{\epsilon, \varepsilon, \tilde{h}^{\epsilon}}_{t} )\dif t+ \frac{1}{\sqrt{\epsilon\varepsilon}}\sigma_{1}(X^{\epsilon, \varepsilon, \tilde{h}^{\epsilon}}_{t}, \sL_{X^{\epsilon, \varepsilon}_{t}}, Y^{\epsilon, \varepsilon, \tilde{h}^{\epsilon}}_{t} ){\cal P}_{1}\dot{h}^{\epsilon}(t)\dif t\\
 &\quad \quad \quad\quad \quad \quad\quad \quad \quad + \frac{1}{\sqrt{\epsilon\varepsilon}}\sigma_{2}(X^{\epsilon, \varepsilon, \tilde{h}^{\epsilon}}_{t}, \sL_{X^{\epsilon, \varepsilon}_{t}}, Y^{\epsilon, \varepsilon, \tilde{h}^{\epsilon}}_{t} ){\cal P}_{2}\dot{h}^{\epsilon}(t)\dif t\\
 &\quad \quad \quad+ \frac{1}{\sqrt{\varepsilon}}\sigma_{1}(X^{\epsilon, \varepsilon, \tilde{h}^{\epsilon}}_{t}, \sL_{X^{\epsilon, \varepsilon}_{t}}, Y^{\epsilon, \varepsilon, \tilde{h}^{\epsilon}}_{t} )\dif W^{1}_{t}
 + \frac{1}{\sqrt{\varepsilon}}\sigma_{2}(X^{\epsilon, \varepsilon, \tilde{h}^{\epsilon}}_{t}, \sL_{X^{\epsilon, \varepsilon}_{t}}, Y^{\epsilon, \varepsilon, \tilde{h}^{\epsilon}}_{t} )\dif W^{2}_{t}, \\
 & \quad \quad \quad\quad \quad \quad\quad \quad \quad \quad \quad \quad\quad \quad \quad\quad \quad \quad \quad \quad \quad\quad \quad \quad\quad \quad \quad \quad \quad \quad\quad  Y^{\epsilon, \varepsilon, \tilde{h}^{\epsilon}}_{0}=y,
 \end{cases}
 \end{align}
where ${\cal P}_{1}: \mR^{d_{1}+d_{2}}\rightarrow \mR^{d_{1}}, {\cal P}_{2}: \mR^{d_{1}+d_{2}}\rightarrow \mR^{d_{2}}$ are two projection operators.
In this part, we will investigate the LDP for Eq. \eqref{31}.   We need some definitions of the theory of LDP.

\bd A nonnegative function $I$ is called a rate function on $C(0,T; \mR^{d})$ if it is lower semicontinuous. Moreover, $I$ is a good rate function  if for each constant $M< \infty,$ the level  set $\{x
\in C(0,T; \mR^{d}): I(x)\leq M\}$ is a compact subset of $C(0,T; \mR^{d}).$
\ed

\bd
Let $I$ be a rate function on $C(0,T; \mR^{d}).$  The family  $\{X^{\epsilon, \varepsilon}:=\Gamma^{\epsilon}(\sqrt{\epsilon}W_{\cdot}, \epsilon^{H}B^{H}_{\cdot}  )\}_{\epsilon> 0}$ of $C(0,T; \mR^{d})-$valued random variables is said to be satisfied a LDP on $C(0,T; \mR^{d})$ with speed $l(\epsilon):=(\epsilon, \epsilon^{2H})$ and rate function $I$ if the following two conditions hold:
\begin{itemize}
\item[$1^{\circ}$] (Upper bound) For each closed subset $F\in C(0,T; \mR^{d}),$
\begin{align}\label{lp1}
\limsup_{\epsilon\rightarrow 0}\epsilon^{2H}\log (P(X^{\epsilon, \varepsilon} \in F)) \leq -\inf_{x\in F}I(x).
\end{align}
\item[$2^{\circ}$] (Lower bound) For each open subset $G\in C(0,T; \mR^{d}),$
\begin{align}\label{lp2}
\liminf_{\epsilon\rightarrow 0}\epsilon\log (P(X^{\epsilon, \varepsilon} \in G))\geq -\inf_{x\in G}I(x).
\end{align}
\end{itemize}
\ed

\bd
Let $I$ be a rate function on $C(0,T; \mR^{d}).$  $\{X^{\epsilon, \varepsilon}:=\Gamma^{\epsilon}(\sqrt{\epsilon}W_{\cdot}, \epsilon^{H}B^{H}_{\cdot}  )\}_{\epsilon> 0}$ is said to be satisfied the Laplace principle upper bound (respectively, lower bound ) on $C(0,T; \mR^{d})$  with speed $l(\epsilon):=(\epsilon, \epsilon^{2H})$ and rate function $I$ if for any bounded continuous function $\rho $ on $C(0,T; \mR^{d}),$
\begin{itemize}
\item[$1^{\circ}$]
\begin{align}\label{k}
\limsup_{\epsilon\rightarrow 0}-\epsilon^{2H}\log \mE\bigg[\bigg.\exp\bigg(\bigg. -\frac{\rho(X^{\epsilon, \varepsilon})}{\epsilon^{2H}} \bigg)\bigg.   \bigg]\bigg. \leq \inf_{x\in C(0,T; \mR^{d})}\{\rho(x)+I(x)\}.
\end{align}
\item[$2^{\circ}$] (Lower bound) For any open subset $G_{i}\in C(0,T; \mR^{d}),i=1,2,$
\begin{align}\label{kk}
\liminf_{\epsilon\rightarrow 0}-\epsilon\log \mE\bigg[\bigg.\exp\bigg(\bigg. -\frac{\rho(X^{\epsilon, \varepsilon})}{\epsilon} \bigg)\bigg.   \bigg]\bigg. \geq \inf_{x\in C(0,T; \mR^{d})}\{\rho(x)+I(x)\}.
\end{align}
\end{itemize}
\ed

\begin{rem}
In these definitions, we use different speeds for the upper bound  and lower  bound, i.e.  the speed in upper bound \eqref{lp1} is $\epsilon^{2H}$, while the speed in lower bound \eqref{lp2} is $\epsilon.$ In fact, this is reasonable. Noting that $H\ge \frac12$  and $\epsilon\in (0,1)$, we therefore have
\begin{align*}
&\limsup_{\epsilon\rightarrow 0}\epsilon\log (P(X^{\epsilon, \varepsilon} \in F))\le \limsup_{\epsilon\rightarrow 0}\epsilon^{2H}\log (P(X^{\epsilon, \varepsilon} \in F)),\\
&\liminf_{\epsilon\rightarrow 0}\epsilon\log (P(X^{\epsilon, \varepsilon} \in G))\le
\liminf_{\epsilon\rightarrow 0}\epsilon^{2H}\log (P(X^{\epsilon, \varepsilon} \in G)).
\end{align*}
This implies the consistency the definition of the LDP of SDEs driven by Brownian motion or fraction Brownian motion, respectively.
\end{rem}

We give the following sufficient condition for  LDP criteria, which  is a version of \cite[Theorem 4.2]{BDM} on the product space.
\bl\label{0101}
  Assume the following conditions hold:
\begin{itemize}
\item[$1^{\circ}$] For any $\epsilon > 0,$ $\Gamma^{\epsilon}: C([0,T];\mR^{d_{1}+d_{2}})\times C(0,T; \mR^{d_{1}})\rightarrow C(0,T; \mR^{d})$ is a measurable mapping.
\item[$2^{\circ}$] Let $\Gamma^{0}: C(0,T;\mR^{d_{1}+d_{2}})\times I^{H+\frac{1}{2}}_{0+}(L^{2}(0,T; \mR^{d_{1}}))\rightarrow C(0,T; \mR^{d})$  be a measurable mapping.
\item[$3^{\circ}$] For every $N>0,$ and  any family $\{h_{\epsilon}; \epsilon >0 \}\subset \mA_{N}, $  $\{\bar{h}_{\epsilon}; \epsilon >0 \}\subset \mathcal{A}_{N}$       satisfying that $h_{\epsilon}\Rightarrow h, \bar{h}_{\epsilon}\Rightarrow \bar{h},\epsilon \rightarrow 0,$ then
    $$\Gamma^{\epsilon}\bigg(\bigg.\sqrt{\epsilon}W_{\cdot}+\int^{\cdot}_{0} \dot{h}_{\epsilon}(\cdot)\dif s, \epsilon^{H}B^{H}_{\cdot}+ R_{H}\bar{h}_{\epsilon}(\cdot)\bigg)\bigg.     \Rightarrow \Gamma^{\circ}\bigg(\bigg.\int^{\cdot}_{0} \dot{h}(s)\dif s, R_{H}\bar{h}(\cdot)  \bigg)\bigg., \epsilon\rightarrow 0.$$
\item[$4^{\circ}$] For every $N>0,$ the set $\{\Gamma^{\circ}\bigg(\bigg.\int^{\cdot}_{0} \dot{h}(s)\dif s, R_{H}\bar{h}(\cdot)\bigg)\bigg.; h\in \mS_{N}, \bar{h}\in \mathcal{S}_{N} \}$ is a compact subset of $ C(0,T; \mR^{d}).$
\end{itemize}
Then the family $\{\Gamma^{\epsilon}(\sqrt{\epsilon}W_{\cdot}, \epsilon^{H}B^{H}_{\cdot})\}$ satisfies a  large deviation principle in $C(0,T; \mR^{d})$ with the rate function $I$ given by
\begin{equation}\label{40}
I(g):= \inf_{\{(h, \bar{h})\in (\mH,\mathcal{H}); g=\Gamma^{0}(\int^{\cdot}_{0}\dot{h}(s)\dif s,  R_{H}\bar{h})\}}\bigg\{\bigg.\frac{1}{2}\int^{T}_{0}|\dot{h}(s)|^{2}\dif s +\frac{1}{2}\|\bar{h}\|^{2}_{\mathcal{H}}   \bigg\}\bigg., g\in C(0,T; \mR^{d}),
 \end{equation}
with $\inf \emptyset =\infty$ by convention.

\el
\begin{proof}
The proof is placed in the appendix.

\end{proof}

 The following lemma is equivalent to the above one.
\bl\label{70}
Assume the following conditions hold:
\begin{itemize}
\item[$1^{\circ}$] Let $\{h^{\epsilon}: \epsilon> 0\}\subset \mA_{N}, \{\bar{h}^{\epsilon}: \epsilon> 0\}\subset \mathcal{A}_{N}. $ For any $\varepsilon_{0} > 0,$
\begin{align*}
&\lim_{\epsilon \rightarrow 0}\bigg[\bigg. d\bigg(\bigg. \Gamma^{\epsilon}\bigg(\bigg.\sqrt{\epsilon}W_{\cdot}+\int^{\cdot}_{0} \dot{h}_{\epsilon}(\cdot)\dif s, \epsilon^{H}B^{H}_{\cdot}+ R_{H}\bar{h}_{\epsilon}(\cdot)\bigg)\bigg., \\
&\quad \quad \quad \quad\quad \quad\quad \quad\quad \quad\Gamma^{\circ}\bigg(\bigg.\int^{\cdot}_{0} \dot{h}^{\epsilon}(s)\dif s, R_{H}\bar{h}^{\epsilon}(\cdot)  \bigg)\bigg.\bigg)\bigg. >\varepsilon_{0}  \bigg]\bigg.\\
&=0,
\end{align*}
where $d(\cdot, \cdot)$ stands for the metric on $C(0,T; \mR^{d}).$
\item[$2^{\circ}$] Let $\{h^{\epsilon}\}\subset \mS_{N}, \{\bar{h}^{\epsilon}\}\subset \mathcal{S}_{N}.$ If $h^{\epsilon}$ converges to some element $h$ in $\mS_{N}$ and $\bar{h}^{\epsilon}$ converges to some element $\bar{h}$ in $\mathcal{S}_{N}, n\rightarrow \infty,$  then $$\Gamma^{\circ}\bigg(\bigg.\int^{\cdot}_{0} \dot{h}^{\epsilon}(s)\dif s, R_{H}\bar{h}^{\epsilon}(\cdot)  \bigg)\bigg.\rightarrow \Gamma^{\circ}\bigg(\bigg.\int^{\cdot}_{0} \dot{h}(s)\dif s, R_{H}\bar{h}(\cdot)  \bigg)\bigg.\, \mbox{in}\, C(0,T; \mR^{d}).$$
\end{itemize}
Set $X^{\epsilon}_{\cdot}=\Gamma^{\epsilon}(\sqrt{\epsilon}W_{\cdot}, \epsilon^{H}B^{H}_{\cdot}). $  Then $\{X^{\epsilon}, \epsilon > 0\}$ satisfies the Laplace principle (hence the LDP) on $ C(0,T; \mR^{d})$ with the rate function $I$ given by  \eqref{40}.
\el
\begin{proof}
Since the proof is similar to that in \cite[Theorem 3.2]{MSZ},  we omit it here.
\end{proof}

 The following  skeleton equation  w.r.t. the slow component of stochastic system \eqref{31} will be used later.

\begin{align}\label{43}
\dif \bar{X}^{\tilde{h}}_{t}&=\bar{f}(t,\bar{X}^{\tilde{h}}_{t}, \sL_{\bar{X}^{0}_{t}} )\dif t
 +  g_{1}(t,\bar{X}^{\tilde{h}}_{t}, \sL_{\bar{X}^{0}_{t}}  ) {\cal P}_{1}\dot{h}\dif t
  +l(t, \sL_{\bar{X}^{0}_{t}} )\dif (R_{H}\bar{h})(t),
 \end{align}
where $\tilde{h}=(h, \bar{h} ),  h\in \mA_{N}, \bar{h}\in \mathcal{A}_{N}.$ Since $\bar{f}, g_{1}$ satisfy Lipschitz condition, Eq. \eqref{43}  has a unique solution $\bar{X}^{\tilde{h}}.$ Define a mapping  $\bar{X}^{\tilde{h}}_{\cdot}=\Gamma^{\circ}\bigg(\bigg.\int^{\cdot}_{0} \dot{h}(s)\dif s, R_{H}\bar{h}(\cdot)  \bigg)\bigg..$  Assume  $\bar{X}^{0}$ satisfies the following equation.

\begin{align}\label{44}
\dif \bar{X}^{0}_{t}&=\bar{f}(\bar{X}^{0}_{t}, \sL_{\bar{X}^{0}_{t}} )\dif t.
 \end{align}
Now,  we are in the position to state our main result.
\bt
Assume $(\mathrm{H}1)$ and $\lim_{\epsilon\rightarrow 0}\frac{\varepsilon}{\epsilon}=0.$  If
$$\sup_{y\in \mR^{d}}|\sigma_{1}(x,\mu,y)|\vee \sup_{y\in \mR^{d}}|\sigma_{2}(x,\mu,y)|\leq C(1+|x|+W_{2}(\mu, \delta_{0})), $$
then   $\{X^{\epsilon, \varepsilon}, \epsilon > 0\}$ satisfies  the LDP on $C(0,T;\mR^{d})$ with the rate function $I$ given by  \eqref{40}.
\et

\subsection{Several priori estimates}Choose a step size $\delta\in (0,1),$ and set  $\bar{t}:= \lfloor\frac{t}{\delta}\rfloor\delta.$  We intend to show  several useful estimates.  Using  similar method as used in Lemma \ref{1+} and following from Lemma \ref{5+},   we derive the following estimates.

\bl \label{45}
Assume $(\mathrm{H}1).$  Then it holds that

\begin{itemize}
\item[{\rm (i)}]$\sup_{\epsilon, \varepsilon \in (0,1)}\sup_{t\in [0, T]}\mE[|X^{\epsilon,\varepsilon}_{t}|^{4}]\leq C(T)(1+|x|^{2}+|y|^{2})$,
\item[{\rm (ii)}]$\sup_{\epsilon, \varepsilon \in (0,1)}\sup_{t\in [0, T]}\mE[|Y^{\epsilon,\varepsilon}_{t}|^{4}]\leq C(T)(1+|x|^{2}+|y|^{2}),$
\item[{\rm (iii)}]$\sup_{\epsilon, \varepsilon \in (0,1)}\mE[|X^{\epsilon,\varepsilon}_{t}- X^{\epsilon,\varepsilon}_{\bar{t}}  |^{2}]\leq C(T)\delta (1+ |x|^{2}+|y|^{2}),$
\item[{\rm (iv)}]$\sup_{h\in \mS_{N}  }\sup_{t\in [0, T]}|\bar{X}^{\tilde{h}}_{t}|^{2}\leq C(N,T)(1+|x|^{2}+\sup_{t\in [0, T]}|\bar{X}^{0}_{t}|^{2}),$
\item[{\rm (v)}]$\mE\int^{T}_{0}|\bar{X}^{\tilde{h}}_{t}- \bar{X}^{\tilde{h}}_{\bar{t}}  |^{2}]\leq C(N, T)\delta (1+ |x|^{2}+|y|^{2}),$

Furthermore, if $\lim_{\epsilon\rightarrow 0}\frac{\varepsilon}{\epsilon}=0,$ we  have
\item[{\rm (vi)}] $\mE[\sup_{t\in [0, T]}|Y^{\epsilon,\varepsilon, \tilde{h}^{\epsilon}}_{t}|^{2}]\leq C(T, N)(1+|x|^{2}+|y|^{2}),$
\item[{\rm (vii)}] $\mE\int^{T}_{0}|X^{\epsilon,\varepsilon, \tilde{h}^{\epsilon}}_{t}|^{2}\dif t\leq C(T, N)(1+|x|^{2}+|y|^{2}),$
\item[{\rm (viii)}] $\mE\int^{T}_{0}|X^{\epsilon,\varepsilon, \tilde{h}^{\epsilon}}_{t}- X^{\epsilon,\varepsilon, \tilde{h}^{\epsilon}}_{\bar{t}}  |^{2}]\leq C(N, T)\delta (1+ |x|^{2}+|y|^{2}).$
\end{itemize}
\el

\subsection{The auxiliary process} Choose a step size $\delta\in (0,1)$ and define an auxiliary process  $\bar{Y}^{\epsilon,\varepsilon}_{t}$ with $\bar Y^{\epsilon,\varepsilon}_{0}
%=Y^{\epsilon,\varepsilon}_{0}
=Y^{\epsilon,\varepsilon, \tilde{h}^{\epsilon}}_{0}=y.$ For $t\in [k\delta, (k+1)\delta], k=1,2,3,\cdots,$
\begin{align}\label{48}
 \bar{Y}^{\epsilon,\varepsilon}_{t}&=\bar{Y}^{\epsilon,\varepsilon}_{k\delta}+ \frac{1}{\varepsilon}\int^{t}_{k\delta}b(X^{\epsilon, \varepsilon, \tilde{h}^{\epsilon}}_{\bar{s}}, \sL_{X^{\epsilon, \varepsilon, \tilde{h}^{\epsilon}}_{\bar{s}}}, \bar{Y}^{\epsilon,\varepsilon}_{s} )\dif s + \frac{1}{\sqrt{\varepsilon}}\int^{t}_{k\delta}\sigma_{1}(X^{\epsilon, \varepsilon, \tilde{h}^{\epsilon}}_{\bar{s}}, \sL_{X^{\epsilon, \varepsilon, \tilde{h}^{\epsilon}}_{\bar{s}}}, \bar{Y}^{\epsilon,\varepsilon}_{s} )\dif W^{1}_{s}\no\\
 &\quad +  \frac{1}{\sqrt{\varepsilon}}\int^{t}_{k\delta}\sigma_{2}(X^{\epsilon, \varepsilon, \tilde{h}^{\epsilon}}_{\bar{s}}, \sL_{X^{\epsilon, \varepsilon, \tilde{h}^{\epsilon}}_{\bar{s}}}, \bar{Y}^{\epsilon,\varepsilon}_{s} )\dif W^{2}_{s}.
\end{align}
We will show the following error estimate between the process $Y^{\epsilon, \varepsilon, \tilde{h}^{\epsilon}}$  and $\bar{Y}^{\epsilon,\varepsilon}.$
\bl\label{49}
Assume $(\mathrm{H}1).$  For any $N, T> 0,$  it holds that
\begin{align*}
\sup_{\epsilon, \varepsilon \in (0,1)}\sup_{t\in [0, T]}\mE[|\bar{Y}^{\epsilon,\varepsilon}_{t}|^{2}]\leq C(T)(1+|x|^{2}+|y|^{2}),
\end{align*}
and
\begin{align*}
\mE\int^{T}_{0}|Y^{\epsilon, \varepsilon, \tilde{h}^{\epsilon}}_{t}-\bar{Y}^{\epsilon,\varepsilon}_{t}|^{2}\dif t\leq C(T,N)(1+|x|^{2}+|y|^{2})(\frac{\varepsilon}{\epsilon}+ \delta).
\end{align*}
\el

\begin{proof}
Similar to  arguments as the proof in Lemma \ref{5+}, we can  obtain the first result.  Next,  we prove the second statement.    Note that
\begin{align*}%\label{50}
Y^{\epsilon, \varepsilon, \tilde{h}^{\epsilon}}_{t}- \bar{Y}^{\epsilon,\varepsilon}_{t}&= \frac{1}{\varepsilon}\int^{t}_{0}[b(X^{\epsilon, \varepsilon, \tilde{h}^{\epsilon}}_{s}, \sL_{X^{\epsilon, \varepsilon, \tilde{h}^{\epsilon}}_{s}}, Y^{\epsilon, \varepsilon, \tilde{h}^{\epsilon}}_{s} )-b(X^{\epsilon, \varepsilon, \tilde{h}^{\epsilon}}_{\bar{s}}, \sL_{X^{\epsilon, \varepsilon, \tilde{h}^{\epsilon}}_{\bar{s}}}, \bar{Y}^{\epsilon,\varepsilon}_{t} )]\dif s \no\\
&+ \frac{1}{\sqrt{\varepsilon}}\int^{t}_{0}[\sigma_{1}(X^{\epsilon, \varepsilon, \tilde{h}^{\epsilon}}_{s}, \sL_{X^{\epsilon, \varepsilon}_{t}}, Y^{\epsilon, \varepsilon, \tilde{h}^{\epsilon}}_{t} )-\sigma_{1}(X^{\epsilon, \varepsilon, \tilde{h}^{\epsilon}}_{\bar{s}}, \sL_{X^{\epsilon, \varepsilon, \tilde{h}^{\epsilon}}_{\bar{s}}}, \bar{Y}^{\epsilon,\varepsilon}_{s} )]\dif W^{1}_{s}\no\\
 &\quad +  \frac{1}{\sqrt{\varepsilon}}\int^{t}_{0}[\sigma_{2}(X^{\epsilon, \varepsilon, \tilde{h}^{\epsilon}}_{s}, \sL_{X^{\epsilon, \varepsilon}_{s}}, Y^{\epsilon, \varepsilon, \tilde{h}^{\epsilon}}_{s} )-\sigma_{2}(X^{\epsilon, \varepsilon, \tilde{h}^{\epsilon}}_{\bar{s}}, \sL_{X^{\epsilon, \varepsilon, \tilde{h}^{\epsilon}}_{\bar{s}}}, \bar{Y}^{\epsilon,\varepsilon}_{s} )]\dif W^{2}_{s}\no\\
 &+ \frac{1}{\sqrt{\epsilon\varepsilon}}\int^{t}_{0}\sigma_{1}(X^{\epsilon, \varepsilon, \tilde{h}^{\epsilon}}_{s}, \sL_{X^{\epsilon, \varepsilon}_{s}}, Y^{\epsilon, \varepsilon, \tilde{h}^{\epsilon}}_{s} ){\cal P}_{1}\dot{h}^{\epsilon}_{s}\dif s\no\\
&+ \frac{1}{\sqrt{\epsilon\varepsilon}}\int^{t}_{0}\sigma_{2}(X^{\epsilon, \varepsilon, \tilde{h}^{\epsilon}}_{s}, \sL_{X^{\epsilon, \varepsilon}_{s}}, Y^{\epsilon, \varepsilon, \tilde{h}^{\epsilon}}_{s} ){\cal P}_{2}\dot{h}^{\epsilon}_{s}\dif s.
\end{align*}
By It\^{o}'s formula,  we have
\begin{align*}%\label{50}
&\frac{\dif }{\dif t}\mE[|Y^{\epsilon, \varepsilon, \tilde{h}^{\epsilon}}_{t}- \bar{Y}^{\epsilon,\varepsilon}_{t}|^{2}]\no\\
&\leq \frac{2}{\varepsilon}\mE\langle Y^{\epsilon, \varepsilon, \tilde{h}^{\epsilon}}_{t}- \bar{Y}^{\epsilon,\varepsilon}_{t},   b(X^{\epsilon, \varepsilon, \tilde{h}^{\epsilon}}_{t}, \sL_{X^{\epsilon, \varepsilon, \tilde{h}^{\epsilon}}_{t}}, Y^{\epsilon, \varepsilon, \tilde{h}^{\epsilon}}_{t} )-b(X^{\epsilon, \varepsilon, \tilde{h}^{\epsilon}}_{\bar{t}}, \sL_{X^{\epsilon, \varepsilon, \tilde{h}^{\epsilon}}_{\bar{t}}}, \bar{Y}^{\epsilon,\varepsilon}_{t} )    \rangle \no\\
&+ \frac{2}{\sqrt{\epsilon\varepsilon}}\mE\langle Y^{\epsilon, \varepsilon, \tilde{h}^{\epsilon}}_{t}- \bar{Y}^{\epsilon,\varepsilon}_{t},  \sigma_{1}(X^{\epsilon, \varepsilon, \tilde{h}^{\epsilon}}_{s}, \sL_{X^{\epsilon, \varepsilon}_{s}}, Y^{\epsilon, \varepsilon, \tilde{h}^{\epsilon}}_{s} ){\cal P}_{1}\dot{h}^{\epsilon}    \rangle\no\\
& + \frac{2}{\sqrt{\epsilon\varepsilon}}\mE\langle Y^{\epsilon, \varepsilon, \tilde{h}^{\epsilon}}_{t}- \bar{Y}^{\epsilon,\varepsilon}_{t},  \sigma_{2}(X^{\epsilon, \varepsilon, \tilde{h}^{\epsilon}}_{s}, \sL_{X^{\epsilon, \varepsilon}_{s}}, Y^{\epsilon, \varepsilon, \tilde{h}^{\epsilon}}_{s} ){\cal P}_{2}\dot{h}^{\epsilon}    \rangle\no\\
&+\frac{1}{\varepsilon}\mE|\sigma_{1}(X^{\epsilon, \varepsilon, \tilde{h}^{\epsilon}}_{t}, \sL_{X^{\epsilon, \varepsilon}_{t}}, Y^{\epsilon, \varepsilon, \tilde{h}^{\epsilon}}_{t} )-\sigma_{1}(X^{\epsilon, \varepsilon, \tilde{h}^{\epsilon}}_{\bar{t}}, \sL_{X^{\epsilon, \varepsilon, \tilde{h}^{\epsilon}}_{\bar{t}}}, \bar{Y}^{\epsilon,\varepsilon}_{t} )|^{2}\no\\
&+\frac{1}{\varepsilon}\mE|\sigma_{2}(X^{\epsilon, \varepsilon, \tilde{h}^{\epsilon}}_{t}, \sL_{X^{\epsilon, \varepsilon}_{t}}, Y^{\epsilon, \varepsilon, \tilde{h}^{\epsilon}}_{t} )-\sigma_{2}(X^{\epsilon, \varepsilon, \tilde{h}^{\epsilon}}_{\bar{t}}, \sL_{X^{\epsilon, \varepsilon, \tilde{h}^{\epsilon}}_{\bar{t}}}, \bar{Y}^{\epsilon,\varepsilon}_{t} )|^{2}=:\sum^{5}_{i=1}M_{i}.
\end{align*}
For $M_{1}, M_{4}, M_{5},$  by $(\mathrm{H}1)$,   we have
\begin{align}\label{51}
M_{1}+ M_{4}+ M_{5}&\leq -\frac{\tilde{\alpha}}{\varepsilon}\mE[|Y^{\epsilon, \varepsilon, \tilde{h}^{\epsilon}}_{t}- \bar{Y}^{\epsilon,\varepsilon}_{t}|^{2}]+\frac{C}{\varepsilon}\mE[|X^{\epsilon, \varepsilon, \tilde{h}^{\epsilon}}_{t}- X^{\epsilon, \varepsilon, \tilde{h}^{\epsilon}}_{\bar{t}}|^{2}]\no\\
& \quad+ \frac{C}{\varepsilon}W_{2}^{2}(\sL_{X^{\epsilon, \varepsilon}_{t}}, \sL_{X^{\epsilon, \varepsilon}_{\bar{t}}}),
\end{align}
where $\tilde{\alpha} \in (0, \alpha).$
For $M_{2}, M_{3}$  by Young's inequality,  we derive
\begin{align}\label{52}
M_{2}+ M_{3}&\leq \frac{\tilde{\alpha}_{1}}{\varepsilon}\mE[|Y^{\epsilon, \varepsilon, \tilde{h}^{\epsilon}}_{t}- \bar{Y}^{\epsilon,\varepsilon}_{t}|^{2}]+ \frac{C}{\epsilon}\mE[(1+|X^{\epsilon, \varepsilon, \tilde{h}^{\epsilon}}_{t}|^{2}+\sL_{X^{\epsilon, \varepsilon}}(|\cdot|^{2}))\|{\cal P}_{1}\||\dot{h}^{\epsilon}_{t}|^{2}]\no\\
&\quad + \frac{C}{\epsilon}\mE[(1+|X^{\epsilon, \varepsilon, \tilde{h}^{\epsilon}}_{t}|^{2}+\sL_{X^{\epsilon, \varepsilon}}(|\cdot|^{2}))\|{\cal P}_{2}\||\dot{h}^{\epsilon}_{t}|^{2}].
\end{align}
where $\tilde{\alpha}_{1} \in (0, \tilde{\alpha}).$
By the above calculations,  we have
\begin{align*}
\frac{\dif }{\dif t}\mE[|Y^{\epsilon, \varepsilon, \tilde{h}^{\epsilon}}_{t}- \bar{Y}^{\epsilon,\varepsilon}_{t}|^{2}]&\leq -\frac{\varsigma}{\varepsilon}\mE[|Y^{\epsilon, \varepsilon, \tilde{h}^{\epsilon}}_{t}- \bar{Y}^{\epsilon,\varepsilon}_{t}|^{2}]+\frac{C}{\varepsilon}\mE[|X^{\epsilon, \varepsilon, \tilde{h}^{\epsilon}}_{t}- X^{\epsilon, \varepsilon, \tilde{h}^{\epsilon}}_{\bar{t}}|^{2}]\\
&\quad +\frac{C}{\varepsilon}\mE[|X^{\epsilon, \varepsilon}_{t}- X^{\epsilon, \varepsilon}_{\bar{t}}|^{2}]+ \frac{C}{\epsilon}\mE[(1+|X^{\epsilon, \varepsilon, \tilde{h}^{\epsilon}}_{t}|^{2}+\sL_{X^{\epsilon, \varepsilon}}(|\cdot|^{2}))\|{\cal P}_{1}\||\dot{h}^{\epsilon}_{t}|^{2}]\no\\
&\quad + \frac{C}{\epsilon}\mE[(1+|X^{\epsilon, \varepsilon, \tilde{h}^{\epsilon}}_{t}|^{2}+\sL_{X^{\epsilon, \varepsilon}}(|\cdot|^{2}))\|{\cal P}_{2}\||\dot{h}^{\epsilon}_{t}|^{2}],
\end{align*}
where $\varsigma=\tilde{\alpha}-\tilde{\alpha}_{1}.$   Due to the comparison theorem,  we derive
\begin{align*}
\mE&[|Y^{\epsilon, \varepsilon, \tilde{h}^{\epsilon}}_{t}- \bar{Y}^{\epsilon,\varepsilon}_{t}|^{2}]\\
&\leq \frac{C}{\varepsilon}\int^{t}_{0}e^{-\frac{\varsigma(t-s)}{\epsilon}}\mE[|X^{\epsilon, \varepsilon, \tilde{h}^{\epsilon}}_{s}- X^{\epsilon, \varepsilon, \tilde{h}^{\epsilon}}_{\bar{s}}|^{2}]\dif s+\frac{C}{\varepsilon}\int^{t}_{0}e^{-\frac{\varsigma(t-s)}{\epsilon}}\mE[|X^{\epsilon, \varepsilon}_{s}- X^{\epsilon, \varepsilon}_{\bar{s}}|^{2}]\dif s \\
& \quad +\frac{C}{\epsilon}\int^{t}_{0}e^{-\frac{\varsigma(t-s)}{\epsilon}}\mE[(1+|X^{\epsilon, \varepsilon, \tilde{h}^{\epsilon}}_{s}|^{2}+\sL_{X^{\epsilon, \varepsilon}}(|\cdot|^{2}))\|{\cal P}_{1}\||\dot{h}^{\epsilon}_{s}|^{2}]\dif s\no\\
&\quad + \frac{C}{\epsilon}\int^{t}_{0}e^{-\frac{\varsigma(t-s)}{\epsilon}}\mE[(1+|X^{\epsilon, \varepsilon, \tilde{h}^{\epsilon}}_{s}|^{2}+\sL_{X^{\epsilon, \varepsilon}}(|\cdot|^{2}))\|{\cal P}_{2}\||\dot{h}^{\epsilon}_{s}|^{2}]\dif s.
\end{align*}
From Fubini's theorem, it holds that
\begin{align*}
\mE&\int^{T}_{0}[|Y^{\epsilon, \varepsilon, \tilde{h}^{\epsilon}}_{t}- \bar{Y}^{\epsilon,\varepsilon}_{t}|^{2}]\\
&\leq \frac{C}{\varepsilon}\int^{T}_{0}\mE[|X^{\epsilon, \varepsilon, \tilde{h}^{\epsilon}}_{s}- X^{\epsilon, \varepsilon, \tilde{h}^{\epsilon}}_{\bar{s}}|^{2}]\bigg(\bigg. \int^{T}_{t} e^{-\frac{\varsigma(t-s)}{\epsilon}}\dif t\bigg)\bigg.\dif s\\
& \quad +\frac{C}{\varepsilon}\int^{T}_{0}\mE[|X^{\epsilon, \varepsilon}_{s}- X^{\epsilon, \varepsilon}_{\bar{s}}|^{2}]\bigg(\bigg. \int^{T}_{t} e^{-\frac{\varsigma(t-s)}{\epsilon}}\dif t\bigg)\bigg.\dif s \\
& \quad +\frac{C}{\epsilon}\int^{T}_{0}\mE[(1+|X^{\epsilon, \varepsilon, \tilde{h}^{\epsilon}}_{s}|^{2}+\sL_{X^{\epsilon, \varepsilon}}(|\cdot|^{2}))|\dot{h}^{\epsilon}_{s}|^{2}]\bigg(\bigg. \int^{T}_{t} e^{-\frac{\varsigma(t-s)}{\epsilon}}\dif t\bigg)\bigg.\dif s\no\\
&\quad + \frac{C}{\epsilon}\int^{T}_{0}\mE[(1+|X^{\epsilon, \varepsilon, \tilde{h}^{\epsilon}}_{s}|^{2}+\sL_{X^{\epsilon, \varepsilon}}(|\cdot|^{2}))|\dot{h}^{\epsilon}_{s}|^{2}]\bigg(\bigg. \int^{T}_{t} e^{-\frac{\varsigma(t-s)}{\epsilon}}\dif t\bigg)\bigg.\dif s\\
& \leq C(T,N)(1+|x|^{2}+|y|^{2})(\delta+ \frac{\varepsilon}{\epsilon}),
\end{align*}
as required.\end{proof}

\section{The proof of LDP}
\bl\label{3.1}
Assume $(\mathrm{H}1)$ and $\lim_{\epsilon\rightarrow 0}\frac{\varepsilon}{\epsilon}=0.$  It holds that
\begin{align}
\lim_{\epsilon\rightarrow 0}\sup_{0\leq t\leq T}\mE|X^{\epsilon,\varepsilon}_{t}- \bar{X}^{0}_{t}|^{2}=0.
\end{align}
\el
\begin{proof}
In the same way as in   \cite [Theorem 2.3]{RSX}, we can get the desired result.
\end{proof}

We shall investigate the LDP by using criteria in Lemma \ref{70}.
The criterion $2^{\circ}$ in Lemma \ref{70} will be shown in the following Theorem.
\bt
Assume $(\mathrm{H1})$. Let $\{h^{n}\}\subset  \mS_{N}, \{\bar{h}^{n}\}\subset  \mathcal{S}_{N}$ such that $h^{n}\rightarrow h $ in $\mS_{N}, $ and $\bar{h}^{n}\rightarrow \bar{h} $ in $ \mathcal{S}_{N},$ as $ n \rightarrow \infty, $ respectively.  Then $\Gamma^{\circ}\bigg(\bigg.\int^{\cdot}_{0} \dot{h}^{n}(s)\dif s, R_{H}\bar{h}^{n}(\cdot)  \bigg)\bigg.\rightarrow \Gamma^{\circ}\bigg(\bigg.\int^{\cdot}_{0} \dot{h}(s)\dif s, R_{H}\bar{h}(\cdot)  \bigg)\bigg.$
in $C(0, T; \mR^{d}).$
\et

\begin{proof}
Let $\tilde{h}^{n}=(h^{n}, \bar h^{n})$ and $\bar{X}^{\tilde{h}^{n}}_{\cdot}=\Gamma^{\circ}\bigg(\bigg.\int^{\cdot}_{0} \dot{h}^{n}(s)\dif s, R_{H}\bar{h}^{n}(\cdot)  \bigg)\bigg..$ Then, $\bar{X}^{\tilde{h}^{n}}_{t}$ solves the following equation:
\begin{align}\label{mm}
\dif \bar{X}^{\tilde{h}^{n}}_{t}&=\bar{f}(\bar{X}^{\tilde{h}^{n}}_{t}, \sL_{\bar{X}^{0}_{t}} )\dif t
 +  g_{1}(\bar{X}^{\tilde{h}^{n}}_{t}, \sL_{\bar{X}^{0}_{t}}  ) {\cal P}_{1}\dot{h}^{n}_{t}\dif t
  +l( \sL_{\bar{X}^{0}_{t}} )\dif (R_{H}\bar{h}^{n})(t).
\end{align}
If $h^{n}\rightarrow h $ in $S_{N} $ and $\bar{h}^{n}\rightarrow \bar{h} $ in $\mS_{N},$ as $ n\rightarrow \infty, $ respectively, it suffices to prove that $\bar{X}^{\tilde{h}^{n}}$ converges strongly to $\bar{X}^{\tilde{h}}$ in $C(0, T; \mR^{d})$ with $\tilde{h}=(h, \bar h),$  as $ n\rightarrow \infty,$ which solves
\begin{align}\label{nn}
\dif \bar{X}^{\tilde{h}}_{t}&=\bar{f}(\bar{X}^{\tilde{h}}_{t}, \sL_{\bar{X}^{0}_{t}} )\dif t
 +  g_{1}(\bar{X}^{\tilde{h}}_{t}, \sL_{\bar{X}^{0}_{t}}  ) {\cal P}_{1}\dot{h}\dif t
  +l(\sL_{\bar{X}^{0}_{t}} )\dif (R_{H}\bar{h})(t).
 \end{align}
By \eqref{mm}, we have
\begin{align}\label{27}
 \bar{X}^{\tilde{h}^{n}}_{t}-\bar{X}^{\tilde{h}^{n}}_{s}&= \int^{t}_{s}\bar{f}(\bar{X}^{\tilde{h}^{n}}_{r}, \sL_{\bar{X}^{0}_{r}})\dif r\no\\
 & \quad + \int^{t}_{s}g_{1}(\bar{X}^{\tilde{h}^{n}}_{r}, \sL_{\bar{X}^{0}_{r}}  ) {\cal P}_{1}\dot{h}^{n}_{r}\dif r+\int^{t}_{s}l( \sL_{\bar{X}^{0}_{r}} )\dif (R_{H}\bar{h}^{n})(r).
\end{align}
From $(\mathrm{H}1)$ and Lemma \ref{45}, we have
\begin{align}\label{54}
 \bigg|\bigg. &\int^{t}_{s}\bar{f}(\bar{X}^{\tilde{h}^{n}}_{r}, \sL_{\bar{X}^{0}_{r}}  ) \dif r\bigg|\bigg.
  \leq C(T,N)\int^{t}_{s}(|\bar{X}^{\tilde{h}^{n}}_{r}|+|\bar{X}^{0}_{r}|)\dif r + C(T,N)(t-s)\no\\
 &  \leq C(T,N)(1+ \sup_{t\in [0, T]}|\bar{X}^{0}_{t}|)(t-s),
\end{align}
and
\begin{align}\label{55}
\bigg|\bigg. &\int^{t}_{s}g_{1}(\bar{X}^{\tilde{h}}_{t}, \sL_{\bar{X}^{0}_{t}} ){\cal P}_{1}\dot{h}^{n}_{r}\dif r\bigg|\bigg. \leq\int^{t}_{s}|g_{1}(\bar{X}^{\tilde{h}}_{t}, \sL_{\bar{X}^{0}_{r}} )||{\cal P}_{1}\dot{h}^{n}_{r}|\dif r\no\\
&\leq \bigg(\bigg.  \int^{t}_{s}|g_{1}(\bar{X}^{\tilde{h}}_{t}, \sL_{\bar{X}^{0}_{t}} )|^{2}\dif r    \bigg)\bigg.^{\frac{1}{2}}\bigg(\bigg. \int^{t}_{s}|{\cal P}_{1}\dot{h}^{n}_{r}|^{2} \dif r   \bigg)\bigg.^{\frac{1}{2}}\no\\
 &  \leq C(T, N)(1+ \sup_{t\in [0, T]}|\bar{X}^{0}_{t}|)(t-s)^{\frac{1}{2}}.
\end{align}
In  view of $(\mathrm{H}1)$ and \eqref{+-}, one has
\begin{align}\label{56}
\bigg|\bigg. &\int^{t}_{s}l( \sL_{\bar{X}^{0}_{r}} )\dif (R_{H}\bar{h}^{n})(r)\bigg|\bigg.=\bigg|\bigg.\int^{t}_{s}l(\sL_{\bar{X}^{0}_{r}} )\int^{r}_{0}\frac{\partial K_{H}}{\partial r}(r,u)(K^{*}_{H}h^{n})(u)\dif  u\dif r\bigg|\bigg.\no\\
& \leq C(T)(1+\sup_{r\in [0, T]}|\bar{X}^{0}_{r}|)\int^{t}_{s}\int^{r}_{0}(\frac{r}{u})^{H-\frac{1}{2}}(r-u)^{H-\frac{3}{2}}|(K^{*}_{H}h^{n})(u)|\dif  u\dif r\no\\
& \leq C(T)(1+\sup_{r\in [0, T]}|\bar{X}^{0}_{r}|)\bigg[\bigg. \int^{s}_{0}u^{\frac{1}{2}-H}|(K^{*}_{H}h^{n})(u)|\bigg(\bigg.\int^{t}_{s}r^{H-\frac{1}{2}} (r-u)^{H-\frac{3}{2}}\dif r \bigg)\bigg.\dif u \no\\
& \quad \quad \quad \quad \quad \quad\quad \quad \quad + \int^{t}_{s}u^{\frac{1}{2}-H}|(K^{*}_{H}h^{n})(u)|\bigg(\bigg.\int^{t}_{u}r^{H-\frac{1}{2}} (r-u)^{H-\frac{3}{2}}\dif r \bigg)\bigg.\dif u \bigg]\bigg..
\end{align}
By the relation $\|K^{*}_{H}h^{n}\|_{L^{2}}= \|h^{n}\|_{\mathcal{H}},$  we derive that
\begin{align}\label{577}
&\int^{s}_{0}u^{\frac{1}{2}-H}|(K^{*}_{H}h^{n})(u)|\bigg(\bigg.\int^{t}_{s}r^{H-\frac{1}{2}} (r-u)^{H-\frac{3}{2}}\dif r \bigg)\bigg.\dif u\no\\
& \leq T^{H-\frac{1}{2}}\int^{t}_{s}(r-s)^{H-\frac{3}{2}}\dif r \int^{s}_{0}u^{\frac{1}{2}-H}|(K^{*}_{H}h^{n})(u)|\dif u\no\\
&\leq \frac{T^{\frac{1}{2}}}{(H-\frac{1}{2})\sqrt{2(1-H)}}\bigg(\bigg. \int^{T}_{0}|(K^{*}_{H}h^{n})(u)|^{2}\dif u \bigg)\bigg.^{\frac{1}{2}}(t-s)^{H-\frac{1}{2}}\no\\
& =\frac{T^{\frac{1}{2}}}{(H-\frac{1}{2})\sqrt{2(1-H)}}\|h^{n}\|_{\mathcal{H}}(t-s)^{H-\frac{1}{2}},
\end{align}
and
\begin{align}\label{588}
&\int^{t}_{s}u^{\frac{1}{2}-H}|(K^{*}_{H}h^{n})(u)|\bigg(\bigg.\int^{t}_{u}r^{H-\frac{1}{2}} (r-u)^{H-\frac{3}{2}}\dif r \bigg)\bigg.\dif u \bigg]\bigg.\no\\
&\leq \frac{T^{H-\frac{1}{2}}}{H-\frac{1}{2}}\int^{t}_{s}u^{\frac{1}{2}-H}(t-u)^{H-\frac{1}{2}}|(K^{*}_{H}h^{n})(u)|\dif u\no\\
&\leq\frac{\sqrt{\mathcal{B}(2-2H, 2H)}T^{H-\frac{1}{2}}}{H-\frac{1}{2}}\|h^{n}\|_{\mathcal{H}}\sqrt{t-s}.
\end{align}
By $(3.1)-(3.6)$ and the fact that $\|h^{n}\|_{\mathcal{H}}\leq \sqrt{2M},$   one can see that $\{\bar{X}^{\tilde{h}^{n}}\}$ is equi-continuous and bounded  in $C(0, T; \mR^{d}),$ which implies $\{\bar{X}^{\tilde{h}^{n}}\}$ is relatively compact in $C(0, T; \mR^{d}).$ Thus, there exists a subsequence still denoted by $\{\bar{X}^{\tilde{h^{n}}}\}$ such that $\{\bar{X}^{\tilde{h^{n}}}\}$ converges to some $\bar X\in C(0, T; \mR^{d}).$

Next, it suffices to prove $\bar{X}=\bar{X}^{\tilde{h}}.$
By $(\mathrm{H}1),$  firstly, we have
\begin{align*}
\bigg|\bigg.\int^{t}_{0}\bar{f}(\bar{X}^{\tilde{h}^{n}}_{s}, \sL_{\bar{X}^{0}_{s}} )\dif s - \int^{t}_{0}\bar{f}(\bar{X}_{s}, \sL_{\bar{X}^{0}_{s}} )\dif s\bigg|\bigg. &\leq C\int^{t}_{0}|\bar{X}^{\tilde{h}^{n}}_{s}-\bar{X}_{s}|\dif s\\
&\leq C(T)\sup_{t\in [0, T]}|\bar{X}^{\tilde{h}^{n}}_{s}-\bar{X}_{s}|\rightarrow 0, n\rightarrow \infty.
\end{align*}
Thus, for any $t\in [0, T], $   we have
\begin{align}\label{59}
\lim_{n\rightarrow \infty}\int^{t}_{0}\bar{f}(\bar{X}^{\tilde{h}^{n}}_{s}, \sL_{\bar{X}^{0}_{s}} )\dif s =\int^{t}_{0}\bar{f}(\bar{X}_{s}, \sL_{\bar{X}^{0}_{s}} )\dif s.
\end{align}
Secondly,
\begin{align}\label{60}
&\bigg|\bigg. \int^{t}_{0}g_{1}(\bar{X}^{\tilde{h}^{n}}_{s}, \sL_{\bar{X}^{0}_{s}}  ) {\cal P}_{1}\dot{h}^{n}_{s}\dif s- \int^{t}_{0}g_{1}(\bar{X}^{\tilde{h}}_{s}, \sL_{\bar{X}^{0}_{s}}  ) {\cal P}_{1}\dot{h}_{s}\dif s\bigg|\bigg.\no\\
& \leq \bigg(\bigg.\int^{t}_{0}|g_{1}(\bar{X}^{\tilde{h}}_{s}, \sL_{\bar{X}^{0}_{s}}  )-g_{1}(\bar{X}^{\tilde{h}}_{s}, \sL_{\bar{X}^{0}_{s}}  ) |^{2}\dif s\bigg)\bigg.^{\frac{1}{2}}\bigg(\bigg.\int^{t}_{0}|\dot{h}^{n}_{s}-\dot{h}_{s}|^{2}\dif s\bigg)\bigg.^{\frac{1}{2}}\no\\
& \leq C(T, N)\sup_{t\in [0, T]}|\bar{X}^{\tilde{h}^{n}}_{t}-\bar{X}_{t}|\rightarrow 0, n\rightarrow \infty.
\end{align}
Then, for any $t\in [0, T], $   we have
\begin{align}\label{61}
\lim_{n\rightarrow \infty}\int^{t}_{0}g_{1}(\bar{X}^{\tilde{h}^{n}}_{s}, \sL_{\bar{X}^{0}_{s}}  ) {\cal P}_{1}\dot{h}^{n}_{s}\dif s\dif s =\int^{t}_{0}g_{1}(\bar{X}^{\tilde{h}}_{s}, \sL_{\bar{X}^{0}_{s}}  ) {\cal P}_{1}\dot{h}_{s}\dif s.
\end{align}
Finally, for any $t\in [0, T],$ we intend to prove
\begin{align}\label{62}
\lim_{n\rightarrow \infty}\int^{t}_{0}l( \sL_{\bar{X}^{0}_{s}} )\dif (R_{H}\bar{h}^{n})(s) =\int^{t}_{0}l( \sL_{\bar{X}^{0}_{s}} )\dif (R_{H}\bar{h})(s)
\end{align}
By Fubini's theorem, we have
\begin{align*}%\label{63}
&\int^{t}_{0}l( \sL_{\bar{X}^{0}_{s}} )\dif (R_{H}\bar{h}^{n})(s) - \int^{t}_{0}l( \sL_{\bar{X}^{0}_{s}} )\dif (R_{H}\bar{h})(s)\no\\
& = \int^{t}_{0}l( \sL_{\bar{X}^{0}_{s}} )\dif s\bigg(\bigg. \int^{s}_{0}\frac{\partial K_{H}}{\partial s}(s,r)[ (K^{*}_{H}h^{n})(r)-(K^{*}_{H}h)(r) ]\dif r   \bigg)\bigg.\no\\
&= C(H)\int^{T}_{0}\bigg[\bigg.1_{[0, t]}(r)r^{\frac{1}{2}- H}\bigg(\bigg.\int^{t}_{r}l( \sL_{\bar{X}^{0}_{s}} )s^{H-\frac{1}{2}}(s-r)^{H-\frac{3}{2}}\dif s   \bigg)\bigg.     \bigg]\bigg.[ (K^{*}_{H}h^{n})(r)-(K^{*}_{H}h)(r) ]\dif r.
\end{align*}
For any unit vector $e \in \mR^{d}, t\in [0,T],$ let
$$\rho_{t}(r) := \bigg.1_{[0, t]}(r)r^{\frac{1}{2}- H}\bigg(\bigg.\int^{t}_{r}l( \sL_{\bar{X}^{0}_{s}} )s^{H-\frac{1}{2}}(s-r)^{H-\frac{3}{2}}\dif s   \bigg)\bigg. e, r\in [0, T].   $$
From $(\mathrm{H}1),$  one has
$$ |\rho_{t}(r)| \leq \frac{C(T, H)}{H+\frac{1}{2}}(1+\sup_{s\in [0, T]}|X^{0}_{s}|)r^{\frac{1}{2}- H}.    $$
This implies $\rho_{t}(\cdot)\in L^{2}(0, T; \mR^{d}).$  Combining \eqref{63} and the condition that $\bar{h}^{n}\rightarrow \bar{h}$  in $\mS_{N},$ we derive
\begin{align}\label{65}
\lim_{n\rightarrow \infty}\int^{t}_{0}l( \sL_{\bar{X}^{0}_{s}} )\dif (R_{H}\bar{h}^{n})(s) =\int^{t}_{0}l( \sL_{\bar{X}^{0}_{s}} )\dif (R_{H}\bar{h})(s).
\end{align}
Taking $n\rightarrow \infty$ in  $\eqref{mm},$  one has  $\bar{X}$ solves \eqref{nn}.
By a standard subsequential argument, we can conclude that the full sequence $\{\bar{X}^{\tilde{h}^{n}}\}$ to $\bar{X}^{\tilde{h}}$ in $\mathcal{E},$  which implies $\bar{X}=\bar{X}^{\tilde{h}}.$   The proof is therefore  complete.
\end{proof}

We now prove the criterion $1^{\circ}$ in Lemma \ref{70}.
\bt
Assume $(\mathrm{H1})$.   Then we have
\begin{align*}
&\lim_{\epsilon \rightarrow 0}P\bigg\{\bigg. d\bigg(\bigg. \Gamma^{\epsilon}\bigg(\bigg.\sqrt{\epsilon}W_{\cdot}+\int^{\cdot}_{0} \dot{h}^{\epsilon}(\cdot)\dif s, \epsilon^{H}B^{H}_{\cdot}+ \frac{\epsilon^{H}}{\epsilon^{\frac{1}{2}}}R_{H}\bar{h}^{\epsilon}(\cdot)\bigg)\bigg., \Gamma^{\circ}\bigg(\bigg.\int^{\cdot}_{0} \dot{h}^{\epsilon}(s)\dif s, R_{H}\bar{h}^{\epsilon}(\cdot)  \bigg)\bigg.\bigg)\bigg. >\varepsilon_{0}  \bigg\}\bigg.\\
&=0.
\end{align*}
\et

\begin{proof}
Note that
\begin{align*}
X^{\epsilon, \varepsilon, \tilde{h}^{\epsilon}}_{t}- \bar{X}^{\tilde{h}^{\epsilon}}_{t}&= \int^{t}_{0}[f_{1}(X^{\epsilon, \varepsilon, \tilde{h}^{\epsilon}}_{s}, \sL_{X^{\epsilon, \varepsilon}}, Y^{\epsilon, \varepsilon, \tilde{h}^{\epsilon}}_{s} )-\bar{f}(\bar{X}^{ \tilde{h}^{\epsilon}}_{s}, \sL_{\bar{X}^{0}_{s}})]\dif s \no\\
&+\int^{t}_{0}[g_{1}(X^{\epsilon, \varepsilon, \tilde{h}^{\epsilon}}_{s}, \sL_{X^{\epsilon, \varepsilon}} )-g_{1}(\bar{X}^{ \tilde{h}^{\epsilon}}_{s}, \sL_{\bar{X}^{0}_{s}})]P_{1}h^{\epsilon}(s)\dif s \no\\
&\quad + \sqrt{\epsilon}\int^{t}_{0}g_{1}(X^{\epsilon, \varepsilon, \tilde{h}^{\epsilon}}_{s}, \sL_{X^{\epsilon, \varepsilon}}, Y^{\epsilon, \varepsilon, \tilde{h}^{\epsilon}}_{s} )\dif W^{1}_{s}\no\\
 & + \int^{t}_{0}[l( \sL_{X^{\epsilon, \varepsilon}_{s}})- l( \sL_{\bar{X}^{0}_{s}})]\dif (R_{H}\bar{h}^{\epsilon})(s)\\
 &+\epsilon^{H} \int^{t}_{0}l( \sL_{X^{\epsilon, \varepsilon}_{s}})\dif B^{H}_{s}.
\end{align*}
Then,  we have
\begin{align*}
|X^{\epsilon, \varepsilon, \tilde{h}^{\epsilon}}_{t}- \bar{X}^{\tilde{h}^{\epsilon}}_{t}|^{2}&\leq  \bigg|\bigg.\int^{t}_{0}[f_{1}(X^{\epsilon, \varepsilon, \tilde{h}^{\epsilon}}_{s}, \sL_{X^{\epsilon, \varepsilon}}, Y^{\epsilon, \varepsilon, \tilde{h}^{\epsilon}}_{s} )-\bar{f}(\bar{X}^{ \tilde{h}^{\epsilon}}_{s}, \sL_{\bar{X}^{0}_{s}})]\dif s\bigg|\bigg.^{2} \no\\
&\quad + \epsilon\bigg|\bigg.\int^{t}_{0}g_{1}(X^{\epsilon, \varepsilon, \tilde{h}^{\epsilon}}_{s}, \sL_{X^{\epsilon, \varepsilon}}, Y^{\epsilon, \varepsilon, \tilde{h}^{\epsilon}}_{s} )\dif W^{1}_{s}\bigg|\bigg.^{2}\no\\
 & + \bigg|\bigg.\int^{t}_{0}[l( \sL_{X^{\epsilon, \varepsilon}_{s}})- l( \sL_{\bar{X}^{0}_{s}})]\dif (R_{H}\bar{h}^{\epsilon})(s)\bigg|\bigg.^{2}\\
 &+\epsilon^{2H}\bigg|\bigg. \int^{t}_{0}l( \sL_{X^{\epsilon, \varepsilon}_{s}})\dif B^{H}_{s}\bigg|\bigg.^{2}\no\\
 &+ \bigg|\bigg.\int^{t}_{0}[g_{1}(X^{\epsilon, \varepsilon, \tilde{h}^{\epsilon}}_{s}, \sL_{X^{\epsilon, \varepsilon}} )-g_{1}(\bar{X}^{ \tilde{h}^{\epsilon}}_{s}, \sL_{\bar{X}^{0}_{s}})]P_{1}h^{\epsilon}(s)\dif s\bigg|\bigg.^{2} \no\\
 &=:\sum^{5}_{i=1}K_{i}(t).
\end{align*}
For $K_{2}(t),$ applying BDG's inequality and Lemma \ref{45},  we have
\begin{align}\label{54+}
 \mE\sup_{0\leq s\leq t}K_{2}(s)\leq C(T)\epsilon(1+|x|^{2}+|y|^{2}).
\end{align}
For $K_{3}(t),$  using  $(\mathrm{H1}),$ \eqref{+-}, the fact that $K^{*}_{H}$ is an isometry between $\mathcal{H}$  and $L^{2}(0, T; \mR^{d_{1}}),$  and the same way  as used in \eqref{577}, \eqref{588}, we arrive at
\begin{align}\label{55Y}
 \mE[\sup_{0\leq s\leq t}K_{3}(s)]\leq C(T, N)\sup_{0\leq s\leq t}\mE[|X^{\epsilon, \varepsilon}_{s}- \bar{X}^{0}_{s}|^{2}].
\end{align}
For $K_{4}(t),$ by $(\mathrm{H1}),$  it yields
\begin{align}\label{56Y}
 \mE[\sup_{0\leq s\leq t}K_{4}(s)]\leq C(T)\epsilon^{2H}(1+|x|^{2}+|y|^{2}).
\end{align}
For $K_{5}(t),$  $(\mathrm{H1})$ implies
\begin{align}\label{566}
\mE\bigg|\bigg.&\int^{t}_{0}[g_{1}(X^{\epsilon, \varepsilon, \tilde{h}^{\epsilon}}_{s}, \sL_{X^{\epsilon, \varepsilon}} )-g_{1}(\bar{X}^{ \tilde{h}^{\epsilon}}_{s}, \sL_{\bar{X}^{0}_{s}})]{\cal P}_{1}h^{\epsilon}(s)\dif s\bigg|\bigg.^{2}\no\\
&\leq  C(T,N)\mE\int^{t}_{0}|X^{\epsilon, \varepsilon, \tilde{h}^{\epsilon}}_{s}-\bar{X}^{ \tilde{h}^{\epsilon}}_{s}|^{2}\dif s  +C(N,T)\mE\int^{t}_{0}|X^{\epsilon, \varepsilon}_{s}-\bar{X}^{0}_{s}|^{2}\dif s.
\end{align}
Next,  we intend to estimate $K_{1}(t).$   By  Lemma \ref{49},  one can derive that
\begin{align}\label{57}
&\mE[\sup_{0\leq s\leq t}K_{1}(s)] \leq C(T)\mE\int^{t}_{0}|f_{1}(X^{\epsilon, \varepsilon, \tilde{h}^{\epsilon}}_{s}, \sL_{X^{\epsilon, \varepsilon}_{s}}, Y^{\epsilon, \varepsilon, \tilde{h}^{\epsilon}}_{s} )-f_{1}(X^{\epsilon, \varepsilon,\tilde{h}^{\epsilon}}_{\bar{s}}, \sL_{X^{\epsilon, \varepsilon}_{\bar{s}}}, \bar{Y}^{\epsilon, \varepsilon}_{s} ) |^{2}\dif s\no\\
&\quad +C(T)\bigg|\bigg.\int^{t}_{0}[f_{1}(X^{\epsilon, \varepsilon, \tilde{h}^{\epsilon}}_{\bar{s}}, \sL_{X^{\epsilon, \varepsilon}_{\bar{s}}}, \bar{Y}^{\epsilon, \varepsilon}_{s}  )- \bar{f}(X^{\epsilon, \varepsilon, \tilde{h}^{\epsilon}}_{\bar{s}}, \sL_{X^{\epsilon, \varepsilon}_{\bar{s}}})]\dif s\bigg|\bigg.^{2}\no\\
&\quad +C(T)\mE\int^{t}_{0}| \bar{f}(X^{\epsilon, \varepsilon, \tilde{h}^{\epsilon}}_{\bar{s}}, \sL_{X^{\epsilon, \varepsilon}_{\bar{s}}})- \bar{f}(X^{\epsilon, \varepsilon, \tilde{h}^{\epsilon}}_{s}, \sL_{X^{\epsilon, \varepsilon}_{s}})|^{2}\dif s\no\\
&\quad +C(T)\mE\int^{t}_{0}|\bar{f}(X^{\epsilon, \varepsilon, \tilde{h}^{\epsilon}}_{s}, \sL_{X^{\epsilon, \varepsilon}_{s}})-\bar{f}(\bar{X}^{ \tilde{h}^{\epsilon}}_{s}, \sL_{\bar{X}^{0}_{s}}) |^{2}\dif s\no\\
& \leq C(T,N)(1+|x|^{2}+|y|^{2})(\delta +\frac{\varepsilon}{\epsilon})  +  C(T,N)\mE\int^{t}_{0}|X^{\epsilon, \varepsilon, \tilde{h}^{\epsilon}}_{s}-\bar{X}^{ \tilde{h}^{\epsilon}}_{s}|^{2}\dif s\no\\
& \quad +C(T)\mE\bigg|\bigg.\int^{t}_{0}[f_{1}(X^{\epsilon, \varepsilon, \tilde{h}^{\epsilon}}_{\bar{s}}, \sL_{X^{\epsilon, \varepsilon}_{\bar{s}}}, \bar{Y}^{\epsilon, \varepsilon}_{\bar{s}}  )- \bar{f}(X^{\epsilon, \varepsilon, \tilde{h}^{\epsilon}}_{\bar{s}}, \sL_{X^{\epsilon, \varepsilon}_{\bar{s}}})]\dif s\bigg|\bigg.^{2}\no\\
& \quad +  C(N,T)\mE\int^{t}_{0}|X^{\epsilon, \varepsilon}_{s}-\bar{X}^{0}_{s}|^{2}\dif s.
\end{align}
Now, we need to estimate  $C(T)\int^{t}_{0}|f_{1}(X^{\epsilon, \varepsilon, \tilde{h}^{\epsilon}}_{\bar{s}}, \sL_{X^{\epsilon, \varepsilon}_{\bar{s}}}, \bar{Y}^{\epsilon, \varepsilon}_{s}  )- \bar{f}(X^{\epsilon, \varepsilon, \tilde{h}^{\epsilon}}_{\bar{s}}, \sL_{X^{\epsilon, \varepsilon}_{\bar{s}}})|^{2}\dif s.$

\begin{align}\label{58}
C(T)&\bigg|\bigg.\int^{t}_{0}[f_{1}(X^{\epsilon, \varepsilon, \tilde{h}^{\epsilon}}_{\bar{s}}, \sL_{X^{\epsilon, \varepsilon}_{\bar{s}}}, \bar{Y}^{\epsilon, \varepsilon}_{s}  )- \bar{f}(X^{\epsilon, \varepsilon, \tilde{h}^{\epsilon}}_{\bar{s}}, \sL_{X^{\epsilon, \varepsilon}_{\bar{s}}})]\dif s\bigg|\bigg.^{2}\no\\
& \leq \frac{C(T)}{\delta}\mE \sum^{k=\lfloor t/\delta\rfloor -1}_{k=0} \bigg|\bigg.\int^{(k+1)\delta}_{k\delta}[f_{1}(X^{\epsilon, \varepsilon, \tilde{h}^{\epsilon}}_{k\delta}, \sL_{X^{\epsilon, \varepsilon}_{k\delta}}, \bar{Y}^{\epsilon, \varepsilon}_{s}  )- \bar{f}(X^{\epsilon, \varepsilon, \tilde{h}^{\epsilon}}_{k\delta}, \sL_{X^{\epsilon, \varepsilon}_{k\delta}})]\dif s\bigg|\bigg.^{2}\no\\
& \quad + C(T)\mE\bigg|\bigg.\int^{t}_{\bar{t}}[f_{1}(X^{\epsilon, \varepsilon, \tilde{h}^{\epsilon}}_{\bar{s}}, \sL_{X^{\epsilon, \varepsilon}_{\bar{s}}}, \bar{Y}^{\epsilon, \varepsilon}_{s}  )- \bar{f}(X^{\epsilon, \varepsilon, \tilde{h}^{\epsilon}}_{\bar{s}}, \sL_{X^{\epsilon, \varepsilon}_{\bar{s}}})]\dif s\bigg|\bigg.^{2}\no\\
&\leq \frac{C(T)\varepsilon^{2}}{\delta^{2}}\max_{0\leq k\leq \lfloor T/\delta\rfloor -1}\bigg\{\bigg.\int^{\frac{\delta}{\varepsilon}}_{0}\int^{\frac{\delta}{\varepsilon}}_{r}\Upsilon_{k}(s,r)\dif s\dif r         \bigg\}\bigg. +C(T,N)\delta(1+|x|^{2}+|y|^{2}),
\end{align}
where for any $0\leq r\leq s\leq\frac{\delta}{\varepsilon},$
\begin{align}\label{59Y}
\Upsilon_{k}(s,r)& := \mE[\langle f_{1}(X^{\epsilon, \varepsilon, \tilde{h}^{\epsilon}}_{k\delta}, \sL_{X^{\epsilon, \varepsilon}_{k\delta}}, \bar{Y}^{\epsilon,\varepsilon}_{\varepsilon s+ k\delta})- \bar{f}(X^{\epsilon, \varepsilon, \tilde{h}^{\epsilon}}_{k\delta}, \sL_{X^{\epsilon, \varepsilon}_{k\delta}}),\no \\
& \quad  f_{1}(X^{\epsilon, \varepsilon, \tilde{h}^{\epsilon}}_{k\delta}, \sL_{X^{\epsilon,\varepsilon}_{k\delta}}, \bar{Y}^{\epsilon,\varepsilon}_{\varepsilon r+ k\delta})- \bar{f}(X^{\epsilon, \varepsilon, \tilde{h}^{\epsilon}}_{k\delta}, \sL_{X^{\epsilon,\varepsilon}_{k\delta}})   \rangle].
\end{align}
For any $s> 0, \mu\in \sP_{2}(\mR^{d}),  x, y\in \mR^{d}, $ consider the following equation:
\begin{align}\label{59+}
 \tilde{Y}^{\varepsilon,s,x,\mu,y}_{t}&=y+ \frac{1}{\varepsilon}\int^{t}_{s} b(x, \mu, \tilde{Y}^{\varepsilon,s,x,\mu,y}_{r} )\dif r  +  \frac{1}{\sqrt{\varepsilon}}\int^{t}_{s} \sigma_{1}(x, \mu, \tilde{Y}^{\varepsilon,s,x,\mu,y}_{r} )\dif W^{1}_{r}\no\\
 &\quad   + \frac{1}{\sqrt{\varepsilon}}\int^{t}_{s} \sigma_{2}(x, \mu, \tilde{Y}^{\varepsilon,s,x,\mu,y}_{r} )\dif W^{2}_{r}, ~~ t\geq s.
\end{align}
Then, we have
\begin{align*}
 \bar{Y}^{\epsilon, \varepsilon}_{t}= \tilde{Y}^{\varepsilon,k\delta,X^{\epsilon, \varepsilon, \tilde{h}^{\epsilon}}_{k\delta}, \sL_{X^{\epsilon, \varepsilon, \tilde{h}^{\epsilon}}_{k\delta}},\bar{Y}^{\epsilon, \varepsilon}_{k\delta}}_{t}, ~~t\in [k\delta, (k+1)\delta].
\end{align*}
It follows from the definition of $\Upsilon_{k}(s,r)$ that
\begin{align*}
\Upsilon_{k}(s,r)& = \mE[\langle f_{1}(X^{\epsilon, \varepsilon, \tilde{h}^{\epsilon}}_{k\delta}, \sL_{X^{\epsilon, \varepsilon}_{k\delta}}, \tilde{Y}^{\varepsilon,k\delta,X^{\epsilon, \varepsilon, \tilde{h}^{\epsilon}}_{k\delta}, \sL_{X^{\epsilon, \varepsilon}_{k\delta},\tilde{h}^{\epsilon}},\bar{Y}^{\epsilon, \varepsilon}_{k\delta}}_{\varepsilon s+ k\delta} )- \bar{f}(X^{\epsilon, \varepsilon, \tilde{h}^{\epsilon}}_{k\delta}, \sL_{X^{\epsilon, \varepsilon}_{k\delta}}), \\
& \quad  f_{1}(X^{\epsilon, \varepsilon, \tilde{h}^{\epsilon}}_{k\delta}, \sL_{X^{\epsilon, \varepsilon}_{k\delta}}, \tilde{Y}^{\epsilon,k\delta, X^{\epsilon,\varepsilon, \tilde{h}^{\epsilon}}_{k\delta}, \sL_{X^{\epsilon,\varepsilon, \tilde{h}^{\epsilon}}_{k\delta}}, \bar{Y}^{\epsilon,\varepsilon}_{k\delta}}_{\varepsilon r+ k\delta} )- \bar{f}(X^{\epsilon, \varepsilon, \tilde{h}^{\epsilon}}_{k\delta}, \sL_{X^{\epsilon, \varepsilon}_{k\delta}})   \rangle].
\end{align*}
Since $\tilde{Y}^{\epsilon,k\delta, x, \mu, y}$ is independent of $\sF_{k\delta},$  and $X^{\epsilon, \varepsilon, \tilde{h}^{\epsilon}}_{k\delta}, \bar{Y}^{\epsilon, \varepsilon}_{k\delta}$ are $\sF_{k\delta}-$measurable,  we have
\begin{align*}
&\Upsilon_{k}(s,r) = \mE\bigg\{\bigg.\mE\bigg[\bigg.\langle f_{1}(X^{\epsilon, \varepsilon, \tilde{h}^{\epsilon}}_{k\delta}, \sL_{X^{\epsilon, \varepsilon}_{k\delta}}, \tilde{Y}^{\varepsilon,k\delta,X^{\epsilon, \varepsilon, \tilde{h}^{\epsilon}}_{k\delta}, \sL_{X^{\epsilon, \varepsilon, \tilde{h}^{\epsilon}}_{k\delta}},\bar{Y}^{\epsilon, \varepsilon}_{k\delta}}_{\varepsilon s+ k\delta} )- \bar{f}(X^{\epsilon, \varepsilon, \tilde{h}^{\epsilon}}_{k\delta}, \sL_{X^{\epsilon, \varepsilon}_{k\delta}}) , \\
&  \quad f_{1}(X^{\epsilon, \varepsilon, \tilde{h}^{\epsilon}}_{k\delta}, \sL_{X^{\epsilon, \varepsilon}_{k\delta}}, \tilde{Y}^{\varepsilon,k\delta,X^{\epsilon, \varepsilon, \tilde{h}^{\epsilon}}_{k\delta}, \sL_{X^{\epsilon, \varepsilon, \tilde{h}^{\epsilon}}_{k\delta}},\bar{Y}^{\epsilon, \varepsilon}_{k\delta}}_{\varepsilon r+ k\delta} )- \bar{f}(X^{\epsilon, \varepsilon, \tilde{h}^{\epsilon}}_{k\delta}, \sL_{X^{\epsilon, \varepsilon}_{k\delta}})    \rangle\bigg|\bigg. \sF_{k\delta}\bigg]\bigg.(\omega)\bigg\}\bigg.\\
& = \mE\bigg\{\bigg.\mE\bigg[\bigg.\langle  f_{1}(X^{\epsilon, \varepsilon, \tilde{h}^{\epsilon}}_{k\delta}(\omega), \sL_{X^{\epsilon, \varepsilon}_{k\delta}}, \tilde{Y}^{\varepsilon,k\delta,X^{\epsilon, \varepsilon, \tilde{h}^{\epsilon}}_{k\delta}(\omega), \sL_{X^{\epsilon, \varepsilon, \tilde{h}^{\epsilon}}_{k\delta}},\bar{Y}^{\epsilon, \varepsilon}_{k\delta}(\omega)}_{\varepsilon s+ k\delta} )- \bar{f}(X^{\epsilon, \varepsilon, \tilde{h}^{\epsilon}}_{k\delta}(\omega), \sL_{X^{\epsilon, \varepsilon}_{k\delta}}) , \\
&  \quad f_{1}(X^{\epsilon, \varepsilon, \tilde{h}^{\epsilon}}_{k\delta}(\omega), \sL_{X^{\epsilon, \varepsilon}_{k\delta}}, \tilde{Y}^{\varepsilon,k\delta,X^{\epsilon, \varepsilon, \tilde{h}^{\epsilon}}_{k\delta}(\omega), \sL_{X^{\epsilon, \varepsilon, \tilde{h}^{\epsilon}}_{k\delta}},\bar{Y}^{\epsilon, \varepsilon}_{k\delta}(\omega)}_{\varepsilon r+ k\delta} )- \bar{f}(X^{\epsilon, \varepsilon, \tilde{h}^{\epsilon}}_{k\delta}(\omega), \sL_{X^{\epsilon, \varepsilon}_{k\delta}})      \rangle\bigg]\bigg.\bigg\}\bigg..
\end{align*}
From the definition of $\tilde{Y}^{\epsilon,s, x, \mu, y},$   we know that
\begin{align}\label{59++}
 \tilde{Y}^{\epsilon,k\delta, x, \mu, y}_{\varepsilon s+ k\delta} &= y+ \int^{s}_{0}b( x,\mu,\tilde{Y}^{\epsilon,k\delta, x, \mu, y}_{\varepsilon r+k\delta} )\dif r + \int^{s}_{0}\sigma_{1}( x,\mu,\tilde{Y}^{\epsilon,k\delta, x, \mu, y}_{\varepsilon r+k\delta} )\dif W^{1,k\delta}_{r}\no\\
 &\quad + \int^{s}_{0}\sigma_{2}( x,\mu,\tilde{Y}^{\epsilon,k\delta, x, \mu, y}_{\varepsilon r+k\delta} )\dif W^{2,k\delta}_{r},
\end{align}
where $W^{1,k\delta}_{r}= \frac{1}{\sqrt{\varepsilon}}(W^{1}_{\varepsilon r+k\delta}- W^{1}_{k\delta}), W^{2,k\delta}_{r}= \frac{1}{\sqrt{\varepsilon}}(W^{2}_{\varepsilon r+k\delta}- W^{2}_{k\delta}).$    Note that
\begin{align}\label{59+++}
 Y^{ x, \mu, y}_{s} &= y+ \int^{s}_{0}b( x,\mu, Y^{ x, \mu, y}_{r} )\dif r + \int^{s}_{0}\sigma_{1}( x,\mu,Y^{ x, \mu, y}_{r} )\dif \tilde{W}^{1}_{r}\no\\
 &\quad + \int^{s}_{0}\sigma_{2}( x,\mu,Y^{ x, \mu, y}_{r} )\dif \tilde{W}^{2}_{r},
\end{align}
By the uniqueness of the solutions of \eqref{59++} and \eqref{59+++}, we have that  $\{\tilde{Y}^{\epsilon,k\delta, x, \mu, y}_{\varepsilon s+ k\delta} )  \}_{0\leq s\leq \frac{\delta}{\varepsilon}}$ and $\{Y^{ x, \mu, y}_{s}  \}_{0\leq s\leq \frac{\delta}{\varepsilon}}$ have the same distribution.  By \cite[Proposition  3.7]{RSX},  it holds that
\begin{align}\label{60Y}
&\Upsilon_{k}(s,r)\no \\
&= \mE\bigg\{\bigg.\tilde{\mE}\bigg[\bigg.\langle  f_{1}(X^{\epsilon, \varepsilon, \tilde{h}^{\epsilon}}_{k\delta}(\omega), \sL_{X^{\epsilon, \varepsilon}_{k\delta}}, Y^{X^{\epsilon, \varepsilon, \tilde{h}^{\epsilon}}_{k\delta}(\omega), \sL_{X^{\epsilon, \varepsilon, \tilde{h}^{\epsilon}}_{k\delta}},\bar{Y}^{\epsilon, \varepsilon}_{k\delta}(\omega)}_{s} )- \bar{f}(X^{\epsilon, \varepsilon, \tilde{h}^{\epsilon}}_{k\delta}(\omega), \sL_{X^{\epsilon, \varepsilon}_{k\delta}}) , \no\\
&  \quad  f_{1}(X^{\epsilon, \varepsilon, \tilde{h}^{\epsilon}}_{k\delta}(\omega), \sL_{X^{\epsilon, \varepsilon}_{k\delta}}, Y^{X^{\epsilon, \varepsilon, \tilde{h}^{\epsilon}}_{k\delta}(\omega), \sL_{X^{\epsilon, \varepsilon, \tilde{h}^{\epsilon}}_{k\delta}},\bar{Y}^{\epsilon, \varepsilon}_{k\delta}(\omega)}_{r} )- \bar{f}(X^{\epsilon, \varepsilon, \tilde{h}^{\epsilon}}_{k\delta}(\omega), \sL_{X^{\epsilon, \varepsilon}_{k\delta}})      \rangle\bigg]\bigg.\bigg\}\bigg.\no\\
& = \mE\bigg\{\bigg.\tilde{\mE}\bigg[\bigg.\tilde{\mE}[\langle f_{1}(X^{\epsilon, \varepsilon, \tilde{h}^{\epsilon}}_{k\delta}(\omega), \sL_{X^{\epsilon, \varepsilon}_{k\delta}}, Y^{X^{\epsilon, \varepsilon, \tilde{h}^{\epsilon}}_{k\delta}(\omega), \sL_{X^{\epsilon, \varepsilon, \tilde{h}^{\epsilon}}_{k\delta}},\bar{Y}^{\epsilon, \varepsilon}_{k\delta}(\omega)}_{s} )|\tilde{\sF}_{r}](\tilde{\omega})\no\\
&\quad \quad \quad\quad \quad \quad\quad \quad \quad\quad \quad \quad\quad \quad \quad\quad \quad \quad\quad \quad \quad\quad \quad \quad - \bar{f}(X^{\epsilon,\varepsilon,\tilde{h}^{\epsilon}}_{k\delta}(\omega), \sL_{X^{\epsilon, \varepsilon}_{k\delta}}), \no\\
&   \quad \quad \quad f_{1}(X^{\epsilon, \varepsilon, \tilde{h}^{\epsilon}}_{k\delta}(\omega), \sL_{X^{\epsilon, \varepsilon}_{k\delta}}, Y^{X^{\epsilon, \varepsilon, \tilde{h}^{\epsilon}}_{k\delta}(\omega), \sL_{X^{\epsilon, \varepsilon, \tilde{h}^{\epsilon}}_{k\delta}},\bar{Y}^{\epsilon, \varepsilon}_{k\delta}(\omega)}_{r} )- \bar{f}(X^{\epsilon,\varepsilon,\tilde{h}^{\epsilon}}_{k\delta}(\omega), \sL_{X^{\epsilon,\varepsilon}_{k\delta}})    \rangle\bigg]\bigg.\bigg\}\bigg.\no\\
&\leq C(T)\mE\bigg\{\bigg.  \tilde{\mE}\bigg[\bigg.1+ |X^{\epsilon,\varepsilon, \tilde{h}^{\epsilon}}_{k\delta}(\omega)|^{2}+ | Y^{X^{\epsilon, \varepsilon, \tilde{h}^{\epsilon}}_{k\delta}(\omega), \sL_{X^{\epsilon, \varepsilon, \tilde{h}^{\epsilon}}_{k\delta}},\bar{Y}^{\epsilon, \varepsilon}_{k\delta}(\omega)}_{r}(\tilde{\omega})|^{2} + \sL_{X^{\epsilon,\varepsilon}_{k\delta}}(|\cdot|^{2})   \bigg]\bigg.e^{-\frac{(s-r)\alpha}{2}}\bigg\}\bigg.\no\\
&\leq C(T)e^{-\frac{(s-r)\alpha}{2}}\mE[1+ |X^{\epsilon,\varepsilon, \tilde{h}^{\epsilon}}_{k \delta}|^{2}+|\hat{Y}^{\varepsilon,\varepsilon}_{k \delta}|^{2}+ \mE(|X^{\epsilon,\varepsilon}_{k \delta}|^{2})+  \mE(|X^{\epsilon,\varepsilon, \tilde{h}^{\epsilon}}_{k \delta}|^{2})]\no\\
&\leq C(T, N)(1+|x|^{2}+|y|^{2})e^{-\frac{(s-r)\alpha}{2}}.
\end{align}
Applying \eqref{58} and \eqref{60}, we conclude
\begin{align}\label{61Y}
C(T)&\bigg|\bigg.\int^{t}_{0}[f_{1}(X^{\epsilon, \varepsilon, \tilde{h}^{\epsilon}}_{\bar{s}}, \sL_{X^{\epsilon, \varepsilon}_{\bar{s}}}, \bar{Y}^{\epsilon, \varepsilon}_{s}  )- \bar{f}(X^{\epsilon, \varepsilon, \tilde{h}^{\epsilon}}_{\bar{s}}, \sL_{X^{\epsilon, \varepsilon}_{\bar{s}}})]\dif s\bigg|\bigg.^{2}\no\\
&\leq C(T,N)(1+|x|^{2}+|y|^{2})(\frac{\varepsilon^{2}}{\delta^{2}}+\frac{\varepsilon}{\delta}+ \delta).
\end{align}
Thus,
\begin{align*}%\label{62Y}
\mE&[\sup_{0\leq s\leq t}|X^{\epsilon, \varepsilon, \tilde{h}^{\epsilon}}_{s}- \bar{X}^{\tilde{h}^{\epsilon}}_{s}|^{2}]\no\\
&\leq C(T, N)(1+|x|^{2}+|y|^{2})(\frac{\varepsilon^{2}}{\delta^{2}}+\frac{\varepsilon}{\delta}+\frac{\varepsilon}{\epsilon}+ \delta+\epsilon)+  C(T,N)\sup_{0\leq t\leq T}\mE[|X^{\epsilon, \varepsilon}_{s}-\bar{X}^{0}_{s}|^{2}]\no\\
&\quad +\int^{t}_{0}\mE[\sup_{0\leq r\leq s}|X^{\epsilon, \varepsilon, \tilde{h}^{\epsilon}}_{r}- \bar{X}^{\tilde{h}^{\epsilon}}_{r}|^{2}]\dif s.
\end{align*}
Taking $\delta= \epsilon^{\frac{2}{3}},$   it yields
\begin{align}\label{63}
\mE&[\sup_{0\leq s\leq t}|X^{\epsilon, \varepsilon, \tilde{h}^{\epsilon}}_{s}- \bar{X}^{\tilde{h}^{\epsilon}}_{s}|^{2}]\no\\
&\leq C(T,N)(1+|x|^{2}+|y|^{2})(\frac{\varepsilon}{\epsilon}+ \varepsilon^{\frac{1}{3}}+\epsilon)+  C(T,N)\sup_{0\leq t\leq T}\mE[|X^{\epsilon, \varepsilon}_{s}-\bar{X}^{0}_{s}|^{2}]\no\\
&\quad +\int^{t}_{0}\mE[\sup_{0\leq r\leq s}|X^{\epsilon, \varepsilon, \tilde{h}^{\epsilon}}_{r}- \bar{X}^{\tilde{h}^{\epsilon}}_{r}|^{2}]\dif s.
\end{align}
By  \eqref{63}, Gronwall's inequality  and Lemma \ref{3.1},   we have for any $\varepsilon_{0}>0$
\begin{align}\label{64}
&P\bigg\{\bigg. d\bigg(\bigg. \Gamma^{\epsilon}\bigg(\bigg.\sqrt{\epsilon}W_{\cdot}+\int^{\cdot}_{0} \dot{h}_{\epsilon}(\cdot)\dif s, \epsilon^{H}B^{H}_{\cdot}+ \frac{\epsilon^{H}}{\epsilon^{\frac{1}{2}}}R_{H}h_{\epsilon}(\cdot)\bigg)\bigg., \Gamma^{\circ}\bigg(\bigg.\int^{\cdot}_{0} \dot{h}^{\epsilon}(s)\dif s, R_{H}\bar{h}^{\epsilon}(\cdot)  \bigg)\bigg.\bigg)\bigg. >\varepsilon_{0}  \bigg\}\bigg.\no\\
&= P\{|X^{\epsilon, \varepsilon, \tilde{h}^{\epsilon}}_{t}- \bar{X}^{\tilde{h}^{\epsilon}}_{t}|> \varepsilon_{0}\}\no\\
&\leq\frac{\mE[\sup_{0\leq t\leq T}|X^{\epsilon, \varepsilon, \tilde{h}^{\epsilon}}_{t}- \bar{X}^{\tilde{h}^{\epsilon}}_{t}|^{2}]}{\varepsilon_{0}}\rightarrow 0, \epsilon\rightarrow 0.
\end{align}

\end{proof}

\section{Appendix}
\subsection{Proof of  Lemma \ref{0101}}
 We now borrow the method in \cite[Theorem 4.4]{BDD} to prove Lemma \ref{0101}. The key is to use the following a variational representation for random functional by making a slight change to that of  \cite[Theorem 3.2]{ZZ}.

\bl\label{a01}
Let $f$ be a bounded Borel measurable function on $\Omega.$ Then it holds that
$$-\log\mE(e^{-f})= \inf_{\tilde{h}=(h, \bar{h}) \in \mathcal{A} \times \mA}\mE(f( \cdot + h(\cdot),\cdot+ R_{H}\bar{h}(\cdot))+\frac{1}{2}\| \bar{h}\|^{2}_{\mathcal{H}}+ \frac{1}{2}\|h\|^{2}_{\mH}).   $$

\el
\begin{proof}[Proof of  Lemma \ref{0101}]
Replacing $f(\cdot)$  by $\frac{\rho\circ \Gamma^{\epsilon}(\sqrt{\epsilon}{\cdot} \epsilon^{H}{\cdot})}{l(\epsilon)} $  in Lemma \ref{a01},  where $\rho$   is real-valued, bounded and continuous function on $\mathcal{E}:=C(0,T; \mR^{d}),$ $l(\epsilon):= \epsilon\, \mbox{or}\, \epsilon^{2H},$ we have

\begin{align}\label{a02}
&-l(\epsilon)\log\mE\bigg[\bigg.e^{-\frac{\rho(X^{\epsilon, \varepsilon})}{l(\epsilon)}}  \bigg]\bigg.= -l(\epsilon)\log\mE\bigg[\bigg.e^{-\frac{\rho\circ \Gamma^{\epsilon}(\sqrt{\epsilon}W_{\cdot}, \epsilon^{H}B^{H}_{\cdot})}{l(\epsilon)}}  \bigg]\bigg.\no\\
&=\inf_{\tilde{h}=(h, \bar{h}) \in \mathcal{A} \times \mA}\mE[\rho\circ \Gamma^{\epsilon}(\sqrt{\epsilon}(W_{\cdot}+h(\cdot)), \epsilon^{H}(B^{H}_{\cdot}+R_{H}\bar{h}(\cdot)))+\frac{1}{2}l(\epsilon)\| \bar{h}\|^{2}_{\mathcal{H}}+ \frac{1}{2}l(\epsilon)\|h\|^{2}_{\mH} ].
\end{align}

The rest of the proof will be divided into two steps.

$\mathbf{Step 1:}$ The upper bound. Without lost of generality,  we assume that $\inf_{x\in \mathcal{E}}\{\rho(x)+I(x) \}<\infty ,$  where $I$ is a rate function given in  \eqref{40}.  Taking $\gamma> 0,$ then there exists $x_{0}\in C(0,T; \mR^{d})$ satisfying
\begin{align}
 \rho(x_{0})+I(x_{0})\leq \inf_{x\in \mathcal{E}}\{\kappa(x)+I(x) \}+ \frac{\gamma}{2}
\end{align}
From \eqref{40}, there exists $(h_{1}, \bar{h}_{1})\in \mathcal{H}\times \mH$ such that $\Gamma^{0}(\int^{\cdot}_{0}\dot{h}_{1}(s)\dif s,  R_{H}\bar{h}
_{1})= x_{0}$ and
$$\frac{1}{2}\|h_{1}\|_{\mH}+ \frac{1}{2}\|\bar{h}_{1}\|_{\mathcal{H}}\leq I(x_{0})+\frac{\gamma}{2}.$$
This together with \eqref{a02} implies

\begin{align*}
&-\epsilon^{2H}\log\mE\bigg[\bigg.e^{-\frac{\rho(X^{\epsilon, \varepsilon})}{\epsilon^{2H}} }  \bigg]\bigg.\\
&\leq\inf_{\tilde{h}=(h, \bar{h}) \in \mathcal{A} \times \mA}\mE[\rho\circ \Gamma^{\epsilon}(\sqrt{\epsilon}(W_{\cdot}+h(\cdot)), \epsilon^{H}(B^{H}_{\cdot}+R_{H}\bar{h}(\cdot)))+\frac{1}{2}\epsilon^{2H}\| \bar{h}\|^{2}_{\mathcal{H}}+ \frac{1}{2}\epsilon\|h\|^{2}_{\mH} ]\no\\
&=\inf_{\tilde{h}=(h, \bar{h}) \in \mathcal{A} \times \mA}\mE[\rho\circ \Gamma^{\epsilon}(\sqrt{\epsilon}(W_{\cdot}+h(\cdot)/\sqrt{\epsilon}), \epsilon^{H}(B^{H}_{\cdot}+R_{H}\bar{h}(\cdot)/\sqrt{\epsilon^{2H}}))\no\\
&\quad \quad \quad\quad\quad \quad \quad\quad\quad \quad \quad\quad\quad \quad \quad\quad\quad\quad \quad \quad\quad+\frac{1}{2}\| \bar{h}\|^{2}_{\mathcal{H}}+ \frac{1}{2}\|h\|^{2}_{\mH} ]\no\\
&\leq \mE[\rho\circ \Gamma^{\epsilon}(\epsilon^{\frac{1}{2}}(W_{\cdot}+h_{1}(\cdot)/\sqrt{\epsilon}), \epsilon^{H}(B^{H}_{\cdot}+R_{H}\bar{h}_{1}(\cdot)/\sqrt{\epsilon^{2H}})) ] +I(x_{0})+\frac{\gamma}{2}.
\end{align*}
By the fact that $\rho$ is bounded and continuous and   taking $\epsilon\rightarrow 0,$  we derive

\begin{align*}
\limsup_{\epsilon\rightarrow 0}-\epsilon^{2H}\log\mE\bigg[\bigg.e^{-\frac{\rho(X^{\epsilon, \varepsilon})}{\epsilon^{2H}}}  \bigg]\bigg.
&\leq \rho\circ \Gamma^{0}(\int^{\cdot}_{0}\dot{h}_{1}(s)\dif s,  R_{H}\bar{h}
_{1}))
+I(x_{0})+\frac{\gamma}{2}\\
&=\rho(x_{0})+ I(x_{0})+\frac{\gamma}{2}\\
&\leq \inf_{x\in \mathcal{E}}\{ \rho(x)+ I(x)\}+\frac{\gamma}{2}.
\end{align*}
Combining  this and the fact that $\gamma$ being arbitrary,  we finish the proof of the upper bound.

$\mathbf{Step 2:}$ The lower bound. Taking  $\gamma> 0,$ by \eqref{a02}, for every $\epsilon > 0,$ there exist $(h^{\epsilon}, \bar{h}^{\epsilon})\in A\times \mA$  such that
\begin{align}\label{a03}
&-\epsilon\log\mE\bigg[\bigg.e^{-\frac{\rho(X^{\epsilon, \varepsilon})}{\epsilon}}  \bigg]\bigg.\no\\
&\geq\inf_{\tilde{h}=(h, \bar{h}) \in \mathcal{A} \times \mA}\mE[\rho\circ \Gamma^{\epsilon}(\sqrt{\epsilon}(W_{\cdot}+h(\cdot)), \epsilon^{H}(B^{H}_{\cdot}+R_{H}\bar{h}(\cdot)))+\frac{1}{2}\epsilon^{2H}\| \bar{h}\|^{2}_{\mathcal{H}}+ \frac{1}{2}\epsilon\|h\|^{2}_{\mH} ]\no\\
& \geq  \mE[\rho\circ \Gamma^{\epsilon}(\epsilon^{\frac{1}{2}}(W_{\cdot}+h^{\epsilon}(\cdot)/\sqrt{\epsilon}), \epsilon^{H}(B^{H}_{\cdot}+R_{H}\bar{h}^{\epsilon}(\cdot)/\sqrt{\epsilon^{2H}}))+\frac{1}{2}\| \bar{h}^{\epsilon}\|^{2}_{\mathcal{H}}+ \frac{1}{2}\|h^{\epsilon}\|^{2}_{\mH} ]-\gamma,
\end{align}
which also implies
\begin{align}\label{a04}
\sup_{\epsilon > 0}\mE\bigg[\bigg.\frac{1}{2}\| \bar{h}^{\epsilon}\|^{2}_{\mathcal{H}}+ \frac{1}{2}\|h^{\epsilon}\|^{2}_{\mH}\bigg]\bigg.\leq 2\|\rho\|_{\infty} +\gamma.
\end{align}
For a given  constant $M>0,$   define the following stopping times.
$$\sigma_{M}^{\epsilon}:=\inf_{t\in [0, T]}\bigg\{\bigg.\frac{1}{2}\| \bar{h}^{\epsilon}1_{[0,t]}\|^{2}_{\mathcal{H}}+ \frac{1}{2}\|1_{[0,t]}h^{\epsilon}\|^{2}_{\mH} \geq M \bigg\}\bigg.\wedge T.  $$
Let $h^{\epsilon, M}(t):=h^{\epsilon}(t)1_{[0, \sigma_{M}^{\epsilon}]}(t), \bar{h}^{\epsilon, M}(t):=\bar{h}^{\epsilon}(t)1_{[0, \sigma_{M}^{\epsilon}]}(t).$  It holds that $h^{\epsilon, M} \in \mA, \bar{h}^{\epsilon, M} \in \mathcal{A}.$  From the Markov inequality and \eqref{a04},  one has
\begin{align}\label{a05}
P(h^{\epsilon}\neq h^{\epsilon, M},h^{\epsilon}\neq h^{\epsilon, M} )\leq P\bigg(\bigg.\frac{1}{2}\| \bar{h}^{\epsilon}\|^{2}_{\mathcal{H}}+ \frac{1}{2}\|h^{\epsilon}\|^{2}_{\mH} \geq M \bigg)\bigg.\leq \frac{2(\|\kappa\|_{\infty}+ \gamma)}{M}.
\end{align}
Moreover, we derive that
\begin{align}\label{a07}
&\rho\circ \Gamma^{\epsilon}(\epsilon^{\frac{1}{2}}(W_{\cdot}+h^{\epsilon}(\cdot)/\sqrt{\epsilon}), \epsilon^{H}(B^{H}_{\cdot}+R_{H}\bar{h}^{\epsilon}(\cdot)/\sqrt{\epsilon^{2H}}))\no\\
&=\rho\circ \Gamma^{\epsilon}(\epsilon^{\frac{1}{2}}(W_{\cdot}+h^{\epsilon, M}(\cdot)/\sqrt{\epsilon}), \epsilon^{H}(B^{H}_{\cdot}+R_{H}\bar{h}^{\epsilon,M}(\cdot)/\sqrt{\epsilon^{2H}}))\no\\
& + \bigg[\bigg.\rho\circ \Gamma^{\epsilon}(\epsilon^{\frac{1}{2}}(W_{\cdot}+h^{\epsilon}(\cdot)/\sqrt{\epsilon}), \epsilon^{H}(B^{H}_{\cdot}+R_{H}\bar{h}^{\epsilon}(\cdot)/\sqrt{\epsilon^{2H}}))\no\\
& \quad \quad \quad \quad - \rho\circ \Gamma^{\epsilon}(\epsilon^{\frac{1}{2}}(W_{\cdot}+h^{\epsilon, M}(\cdot)/\sqrt{\epsilon}), \epsilon^{H}(B^{H}_{\cdot}+R_{H}\bar{h}^{\epsilon,M}(\cdot)/\sqrt{\epsilon^{2H}}))  \bigg]\bigg.1_{\{\bar{h}^{\epsilon}\neq \bar{h}^{\epsilon, M},h^{\epsilon}\neq h^{\epsilon, M}\} }\no\\
&\geq \rho\circ \Gamma^{\epsilon}(\epsilon^{\frac{1}{2}}(W_{\cdot}+h^{\epsilon, M}(\cdot)/\sqrt{\epsilon}), \epsilon^{H}(B^{H}_{\cdot}+R_{H}\bar{h}^{\epsilon}(\cdot)/\sqrt{\epsilon^{2H}}))-2 \|\rho\|_{\infty}1_{\{\bar{h}^{\epsilon}\neq \bar{h}^{\epsilon, M},h^{\epsilon}\neq h^{\epsilon, M}\} },
\end{align}
By \eqref{-+} , one can see that
\begin{align}\label{a08}
&\|\bar{h}^{\epsilon}\|_{\mathcal{H}}^{2}=\|K_{H}^{*}\bar{h}^{\epsilon}\|^{2}_{L^{2}}\geq \|K_{H}^{*}\bar{h}^{\epsilon, M}\|^{2}_{L^{2}}= \|\bar{h}^{\epsilon, M}\|_{\mathcal{H}}^{2},\no\\
&\|h^{\epsilon}\|_{\mathcal{H}}^{2}=\|\dot{h}^{\epsilon}\|^{2}_{L^{2}}\geq \|\dot{h}^{\epsilon, M}\|^{2}_{L^{2}}= \|h^{\epsilon, M}\|_{\mathcal{H}}^{2}.
\end{align}
From  \eqref{a03} \eqref{a05}, \eqref{a07} and \eqref{a08}, we therefore have
\begin{align}\label{a09}
&-\epsilon\log\mE\bigg[\bigg.e^{-\frac{\kappa(X^{\epsilon, \varepsilon})}{\epsilon}}  \bigg]\bigg.\no\\
& \geq  \mE[\rho\circ \Gamma^{\epsilon}(\epsilon^{\frac{1}{2}}(W_{\cdot}+h^{\epsilon}(\cdot)/\sqrt{\epsilon}), \epsilon^{H}(B^{H}_{\cdot}+R_{H}\bar{h}^{\epsilon}(\cdot)/\sqrt{\epsilon^{2H}}))\no\\
&\quad +\frac{1}{2}\| \bar{h}\|^{2}_{\mathcal{H}} +\frac{1}{2}\|h\|^{2}_{\mH} ]-\frac{2\|\rho\|_{\infty}(2\|\kappa\|_{\infty}+\gamma)}{M}-\gamma.
\end{align}
Since $M, \gamma$ are arbitrary, in order to prove the lower bound, it suffices to show
\begin{align}\label{a10}
&\liminf_{\epsilon \rightarrow 0}\mE[\kappa\circ \Gamma^{\epsilon}(\epsilon^{\frac{1}{2}}(W_{\cdot}+h^{\epsilon}(\cdot)/\sqrt{\epsilon}), \epsilon^{H}(B^{H}_{\cdot}+R_{H}\bar{h}^{\epsilon}(\cdot)/\sqrt{\epsilon^{2H}}))+\frac{1}{2}\| \bar{h}\|^{2}_{\mathcal{H}} +\frac{1}{2}\|h\|^{2}_{\mH} ]\no\\
&\geq \inf_{x\in \mathcal{E}}\{\kappa(x)+I(x)  \}.
\end{align}
Since
$$\frac{1}{2}\| \bar{h}^{\epsilon, M}\|^{2}_{\mathcal{H}}\leq M,  \frac{1}{2}\|h^{\epsilon, M}\|^{2}_{\mH} \leq M,     $$
we can extract a (not relabelled) subsequence such that $\bar{h}^{\epsilon, M}$ converges to $\bar{h}$ in distribution and  $h^{\epsilon, M}$ converges to $h$ in distribution.   Then, we obtain
\begin{align*}
&\liminf_{\epsilon \rightarrow 0}\mE[\kappa\circ \Gamma^{\epsilon}(\epsilon^{\frac{1}{2}}(W_{\cdot}+h^{\epsilon}(\cdot)/ \sqrt{\epsilon}, \epsilon^{H}(B^{H}_{\cdot}+R_{H}\bar{h}^{\epsilon}(\cdot)/\sqrt{\epsilon^{2H}}))+\frac{1}{2}\| \bar{h}\|^{2}_{\mathcal{H}} +\frac{1}{2}\|h\|^{2}_{\mH} ]\no\\
&\geq    \mE[\kappa\circ \Gamma^{0}(\dot{h}(\cdot), R_{H}\bar{h}(\cdot))+\frac{1}{2}\| \bar{h}\|^{2}_{\mathcal{H}} +\frac{1}{2}\|h\|^{2}_{\mH} ]   \no\\
&\geq \inf_{(x, h,\bar{h})\in \mathcal{E}
\times \mH\times \mathcal{H}}\mE[\kappa(x)+\frac{1}{2}\| \bar{h}\|^{2}_{\mathcal{H}} +\frac{1}{2}\|h\|^{2}_{\mH} ]\no\\
&\geq \inf_{x\in \mathcal{E}}\{\kappa(x)+I(x)  \}.
\end{align*}
The proof is complete.

\end{proof}

\end{document}